\documentclass{jssc}        

\usepackage{amsmath,amssymb}
\usepackage{graphicx}       
\usepackage{natbib}
\usepackage{algorithm,algorithmic}
\usepackage{xcolor}
\usepackage{verbatim}
\usepackage{multirow}
\usepackage{accents}
\usepackage{enumitem}   

\def\cdd{\mbox{\boldmath$\cdot$}~}

\makeatletter
\def\@oddfoot{\hfill}
\newcount\shumeicount
\def\setshumei#1#2#3{%
  \shumeicount=\count0
  \def\@oddhead{%
    \raise-5pt\hbox to0pt{\vrule width\hsize height 0pt depth 0.4pt\hss}\relax
    \ifnum \shumeicount=\count0
      \raise-7pt\hbox to0pt{\vrule width\hsize height 0pt depth 0.4pt\hss}\relax
      #1
    \else
      \ifodd\count0
        #2
      \else
        #3
       \fi
     \fi
  }%
}
\makeatother
\makeatletter
\def\@oddfoot{\hfill}
\newcount\shujiaocount
\def\setshujiao{%
  \shujiaocount=\count0
  \def\@oddfoot{%
      \ifodd\count0
      \else
      \fi
  }%
}
\makeatother
\def\title#1#2#3#4{{
  \vspace*{0.3cm}
  \begin{flushleft} \Large\bf #1\end{flushleft}
  \vspace*{-0.2cm}
      \begin{flushleft}
      \bf #2
      \end{flushleft}
      \footnotetext{\hspace{-6mm} #3\\ #4}}}

\def\dshm#1#2#3#4
{\setshumei{(#1) #2:
{\thepage--\pageref{LastPage}}\hfill}
            {\hfill {\small #3}\hfill\hbox to0pt{\hss\thepage}}
            {\hbox to0pt{\thepage\hss}\hfill {\small #4}\hfill
            }
            \setshujiao}
\def\drd#1#2
{{\vskip 1cm\small \begin{flushleft}
 #1 \\
 #2\\
\end{flushleft}}}

\def\hat{\widehat}

\def\bar{\overline}
\def\epsilon{\varepsilon}

\begin{document}

\title{Contest for system observability as an infinitely repeated game$^*$}
{\uppercase{Xu} Yueyue \cdd \uppercase{Zhou}
Panpan \cdd \uppercase{Wang} Lin \cdd \uppercase{Liu} Zhixin \cdd \uppercase{Hu} Xiaoming}
{\uppercase{Xu} Yueyue \\
Department of Automation, Shanghai Jiao Tong University, Shanghai 200240, China; Department of Mathematics, KTH Royal Institute of Technology, Stockholm 10044, Sweden.  Email: merryspread99@sjtu.edu.cn.\\
\uppercase{Wang} Lin (Corresponding author) \\
Department of Automation, Shanghai Jiao Tong University, Shanghai 200240, China.  Email: wanglin@sjtu.edu.cn. \\
  \uppercase{Zhou} Panpan \cdd \uppercase{Hu} Xiaoming \\
Department of Mathematics, KTH Royal Institute of Technology, Stockholm 10044, Sweden.  Emails:  panpanz@kth.se; hu@kth.se. \\
\uppercase{Liu} Zhixin \\
Key Laboratory of Systems and Control, Academy of Mathematics and Systems Science, CAS, Beijing 100190, China.  Email: lzx@amss.ac.cn.} 
{$^*$This work was supported in part by the National Natural Science Foundation of China with No. 62373245, in part by the National Key  Research and Development Program of China with No. 2023YFB4706800, and in part by the Dawn Program of Shanghai Education Commission, China. Partial results of this paper were presented at IFAC World Congress 2023.\\} 

\drd{DOI: }{Received: x x 2025}{ / Revised: x x 2025}


\dshm{2025}{XX}{\uppercase{Contest for system observability as an infinitely repeated
game}}{\uppercase{Xu Yueyue}, et al.}

\Abstract{This paper studies a system security problem in the context of observability based on a two-person noncooperative infinitely repeated game. Both the attacker and the defender have means to modify the dimension of the unobservable subspace, which is set as the value function. Utilizing tools from geometric control, we construct the best response sets considering one-step and two-step optimality respectively to maximize or minimize the value function. 
We establish a unified necessary-and-sufficient condition for Nash equilibrium that holds for both one-step and two-step optimizations. Our analysis further uncovers two evolutionary patterns, lock and loop modes, and shows an asymmetry between defense and attack. The defender can lock the game into equilibrium, whereas the attacker can disrupt it by sacrificing short-term utility for longer-term advantage. Six representative numerical examples corroborate the theoretical results and highlight the complexity of possible game outcomes.}      




\section{Introduction}
In recent years, more and more attention has been paid to the security of control systems. Remote sensors are vulnerable to attacks, which intend to deteriorate system performance by manipulating the transmitted data while remaining stealthy \cite{niu2023innovation}. Those attacks may have severe consequences. For example, in the cyber attack on the Ukrainian power grid in 2015, the attacker operated several of the circuit breakers in the grid and jammed the communication network to cause a large scale blackout and keep the operators unaware \cite{soltan2018react}. If the attack can be detected in time, the damage can be reduced. 

In this paper, we want to study attacks against observability. Observability is a critical aspect of system performance. When observability is destroyed, not only the aforementioned attacks are harder to be detected, but there are more severe consequences. For example, observers relying on observability of the system become unusable \cite{bernard2022observer} and operators are unable to accurately recover the true state of the system \cite{showkatbakhsh2020securing}. Plus, systems that lose observability may be more susceptible to various forms of stealthy attacks \cite{liu2024performance}. Researchers have recognized the significance of study on observability under attack in control systems. Regarding state reconstruction after attacks, observability of linear systems has been studied for the case where an attacker modified some outputs \cite{chong2015obser} (along with similar findings by \cite{fawzi2014secure} and \cite{shoukry2015event}).
If more than half of the sensors were attacked, accurately reconstructing the initial state would become impossible. In \cite{mitra2018distributed}, the concept of eigenvalue observability is introduced to estimate locally undetectable states caused by attacks when a single node can exchange information with its neighbours. Later the same authors develop a fully distributed algorithm that successfully reconstructed the system state despite the presence of sensor attacks within the network \cite{mitra2019byzantine}. In \cite{zhang2023observability}, the robustness of observability of a linear time-invariant system under sensor failures is studied from a computational perspective.

However, it should be emphasized that the above research focuses on whether the system can reconstruct the initial state qualitatively. They did not study dimension change of the unobservable subspace quantitatively, which can greatly impact the effectiveness of stealthy attacks. In this context increasing the dimension of the unobservable subspace can facilitate covert attacks, since the stealthy attack vector space is getting larger \cite{kim2014subspace}, \cite{zhao2019sparse}, \cite{maccarone2018uncovering}. In \cite{kim2014subspace} and \cite{zhao2019sparse}, subspace methods are used to construct undetectable attacks, yet without altering the dimension of the unobservable subspace. In \cite{maccarone2018uncovering}, the attacker enhances the effect of stealthy attack by masking sensor signals to increase the dimension of the unobservable subspace. 

In recent years, the game approach has shed new light on the system security problem \cite{zhoupanpan}, \cite{zhang2025}. When the attacker and the defender have limited information, the game becomes partially observable. In \cite{horak2019optimizing}, a partially observable stochastic game is studied and the authors give a representation of uncertainty encountered by the defender. In \cite{zheng2022stackelberg}, an $\epsilon$-Stackelberg partially observable game model is built to prevent state information leakage. However, observability is considered as a game setting in the above research and few game research directly involves confrontation on the observability between attackers and defenders. In \cite{maccarone2020game}, a game approach is used to study the stealthy attack problem in which the attacker masked the sensors and the defender reinforced the sensors. However, in the utility of the players, observability is abstracted into different values without considering the relationship between observability and structural matrices of the system \cite{maccarone2020game}. 

Repeated game is a game model where participants engage in the same basic game for many times. In the system security field, repeated game is a potent tool to analyze the dynamic interplay between attackers and defenders and help both sides to design strategies in their favor. A finite repeated security game is constructed in \cite{nguyen2019deception} to study how an attacker manipulates attack data to mislead the defender, thereby influencing the defender's learning process to make the game more favorable to itself. In \cite{balaji2019design}, an infinitely repeated game and cooperation method are designed to detect malicious nodes and improve energy efficiency. In response to the non-convex infinitely repeated game problem, algorithms are constructed to select both the optimal security strategies that necessitate monitoring and those strategies that do not \cite{aziz2020resilience}.

Additionally, most existing game models use continuous value functions such as quadratic functions, and equilibrium solutions are obtained through methods such as dynamic programming in \cite{xie2025} and Q-learning in \cite{rizvi2025}. However, these methods are not suitable for discrete value functions, such as the dimension of the unobservable subspace.

To overcome the above limitations, this study models the confrontation
between the attacker and the defender as an infinitely repeated game.
The attacker attempts to undermine system observability by maximizing
the dimension of the unobservable subspace, whereas the defender seeks to minimize that dimension so as to preserve observability. The main contributions are summarized below.

(1) We introduce the dimension of the unobservable subspace as the value function and develop closed-form best response algorithms for both players. By combining geometric control computations of controlled invariant subspaces with game model, this work quantifies the adversarial contest over system observability.
In contrast to~\cite{xu2023}, which addresses only one-step
optimality, our framework establishes an extended-horizon optimization formulation and derives analytical expressions.

(2) We establish a unified necessary-and-sufficient condition for Nash equilibrium that holds for both one-step and two-step optimizations, expressed as a concise equality test. Because two players move in alternation, extending the optimization horizon beyond two steps provides no additional strategic information. Previous studies have not provided such a unified equilibrium criterion. 

(3) We discover the defense–attack asymmetry and its resulting game-outcome patterns. The theoretical and experimental analysis shows that the defender can lock the game into equilibrium, whereas the attacker can disrupt it by postponing immediate gains. Two evolutionary paths, lock mode and loop mode, have been identified, providing new theoretical guidance for designing practical security counter-strategies. This finding highlights a fundamental asymmetry neglected in earlier work that usually assumes symmetric influence.

The rest of paper is organized as follows. In Section 2, we derive an observability equivalent system and formulate the game problem. In Section 3, derivations and algorithms are given to get the best response sets for both players. In Section 4, we derive the Nash equilibrium under both one-step and two-step optimality, then refine the game outcome modes and analyze the equilibrium characterization.
In Section 5, examples are given to illustrate various game results. Section 6 is a brief conclusion.
A summary of notations is provided in Table 1.

\begin{table}
\begin{center}
\centerline{\small {\bf Table 1}\ \  Notations}
\vskip0.1cm
\label{tab1}
\begin{tabular}{c c}
\hline
Notations & Definitions \\
\hline
$\mathbb{R}$ & set of real numbers \\
$\mathbb{R}^{n}$ & set of n-dimensional real vectors \\
$\mathbb{R}^{m\times n}$ & set of $m\times n$-dimensional real matrices \\
$\operatorname{Im}A, A\in \mathbb{R}^{m\times n}$ & image space, $\{v\in \mathbb{R}^{m}: v=A q,\forall \, q \in \mathbb{R}^{n}\}$ \\
$\operatorname{Ker}A, A\in \mathbb{R}^{m\times n}$ & kernel space, $\{w\in \mathbb{R}^{n}: A w=0\}$ \\
$\operatorname{pinv}(A)$ & pseudo inverse, $\left(A^\top A\right)^{-1} A^\top$\\
 $\operatorname{Col}_{k}(A)$ & the $k$-th column of matrix $A$ \\
$\delta_{n}^{k}$  & $\operatorname{Col_{k}}\left(I_{n}\right)$, the $k$-th column of $I_{n}$\\
$\oplus$ & direct sum \\
\hline
\end{tabular}
\end{center}
\end{table}

\section{Problem formulation}
\subsection{Modeling of the system}
Consider the following linear system
\begin{equation}\label{1}
\begin{split}
\dot{x}&=Ax+B_du_d+B_au_a,\\
y&=Cx,
\end{split}
\end{equation}
where $x \in \mathbb{R}^{n_0}$, and $y \in \mathbb{R}^{m}$ are state and output of the system respectively; $A \in \mathbb{R}^{n_0\times n_0}$ is system matrix and $C \in \mathbb{R}^{m\times n_0}$ is output matrix. $u_d \in \mathbb{R}^{m}$ and $u_a\in \mathbb{R}^{k}$ are two input channels controlled by the defender and the attacker respectively. 

In this paper, it is assumed that the attacker wants to destroy observability and maximize the unobservable subspace by injecting feedback-data using $u_a$. To the contrary, the system defender wants to protect the system observability and minimize the unobservable subspace via the input $u_d$.
Let us first recall $\mathcal{V}^*$ space and friend matrices.
\begin{definition}\label{defi_1}
$\mathcal{V}$ is a controlled invariant (or $(A,B)$-invariant) subspace if there exists a matrix $F$ such that $(A+BF)\mathcal{V}\subseteq \mathcal{V}$. Such an $F$ is called a friend matrix of $\mathcal{V}$ and we denote the set of friend matrices by $\mathcal{F}(\mathcal{V})$. Define controlled invariant subspaces contained in space $Z$ as $S(Z)$. In $S(\operatorname{Ker}\,C)$, there is a maximal one which is denoted as $\mathcal{V}^*$.
\end{definition}
In this paper, we assume the square system $(A,B_d,C)$ has the relative degree $(r_1, . . . , r_m)$, which reflects the order of differentiation needed in order to have the input $u_d$ explicitly appearing in the output $y$. It is also assumed that the attacker is employing a specific kind of stealthy attack, namely a zero-dynamics attack \cite{Teixeira2012}, which requires $\mbox{Im}\,B_a  \subseteq \mathcal{V}^*$, where $\mathcal{V}^*$ is the maximal $(A,B_d)$-invariant subspace in Ker $C$.

With the relative degree, it is well known that after a coordinate change $x \rightarrow\left[\begin{array}{l}
z \\
\xi
\end{array}\right]$, where $\xi=\left[\xi_{1}^{1},\xi_{2}^{1},\cdots, \xi_{r_{1}}^{1}, \cdots, \xi_{1}^{m}, \xi_{2}^{m},\cdots, \xi_{r_{m}}^{m}\right]^\top$, system \eqref{1} can be rewritten as the normal form
\begin{align*}
\dot{z}&=N z+E \xi+B_a'u_a, \\
\dot{\xi}_{1}^{i} &=\xi_{2}^{i},\\
& \,\vdots \\
\dot{\xi}_{r_{i}-1}^{i} &=\xi_{r_{i}}^{i},\\
\dot\xi_{r_{i}}^{i}&=R_{i} z+S_{i}\xi+c_{i} A^{r_{i}-1} B_{d} u_d, \\
y_{i}&=\xi_{1}^{i},
\end{align*}
where $i=1, \ldots, m , y=[y_1,y_2,\cdots,y_m]^\top$, $B_a'$ is determined by $B_{a}$.
Since the attacker can in essence only change the zero dynamics of the system, the attack is a zero-dynamics attack.

\begin{figure}
\centering \subfigure[Model in the normal form]{
\includegraphics[scale=0.32]{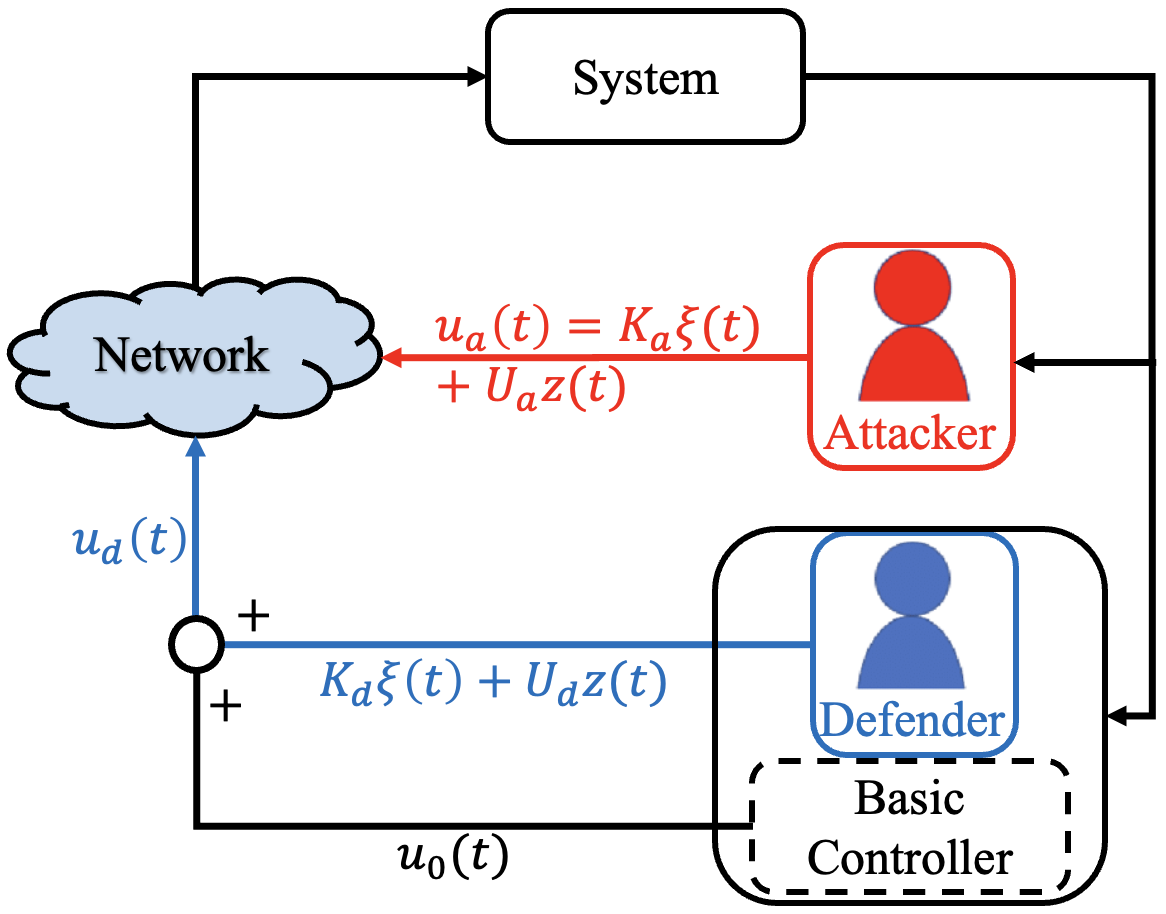}
\label{fig:sub1} } \subfigure[Model in the observability equivalent system]{
\includegraphics[scale=0.4]{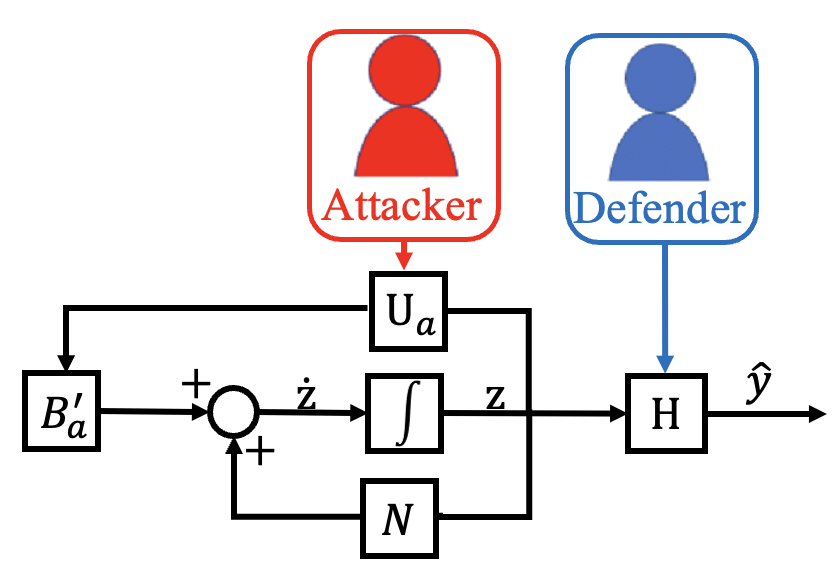}
\label{fig:sub2} }
\centerline {\small \parbox[t]{1.7cm}{\bf Figure 1}
\parbox[t]{10.5cm}
 { Comparison of the model in the normal form (a) and in the observability equivalent system (b).}}
\end{figure}

The complete evolution equation of $\xi_{r_{i}}^{i}$ can be written as
\begin{equation}
[\dot{\xi}_{r_{1}}^{1^\top},\dot{\xi}_{r_{2}}^{2^\top},\cdots,\dot{\xi}_{r_{m-1}}^{m-1^\top},\dot{\xi}_{r_{m}}^{m^\top}]^\top=Rz+ S\xi+L u_d,
\end{equation}
where $R \in \mathbb{R}^{m\times(n_0-s)}$, $S \in \mathbb{R}^{m\times s}$, $L\in \mathbb{R}^{m\times m}$, $s=\sum_{i=1}^{m} r_{i}$.
We define the feedback controls as
\begin{equation}
u_d=K_d\xi+U_dz+u_0,\,u_a=K_a\xi+U_az,
\end{equation}
where $K_d \in \mathbb{R}^{m\times s}$, $U_d \in \mathbb{R}^{m \times (n_0-s)}$ are determined by the defender; $K_a \in \mathbb{R}^{k\times s}$, $U_a \in \mathbb{R}^{k\times (n_0-s)}$ are determined by the attacker; $u_0$ is the input that maintains the normal operation of the system and is controlled by a basic controller, which is illustrated in Figure \ref{fig:sub1}.

Then we have
\begin{equation}
[\dot{\xi}_{r_{1}}^{1^\top},\cdots,\dot{\xi}_{r_{m}}^{m^\top}]^\top=(R+LU_d)z+(S+LK_d)\xi+Lu_0.
\end{equation}
Considering the scenario where $u_0$ does not incorporate state feedback, it consequently does not affect the observability of the system. Thus $u_0$ is omitted in the following derivation of the observability equivalent system. Define 
\begin{equation}
R+LU_d \triangleq {H}.
\end{equation}
Since $L$ is a non-singular matrix according to the definition of relative degree, $\forall {H}\in \mathbb{R}^{m\times n}$, there exists $U_d=L^{-1}\left({H}-R\right)$. Thus the defender can completely control ${H}$. If the system is observable, when $y \equiv \mathbf{0}$, there is $\left[\begin{array}{l}
z \\
\xi
\end{array}\right] \equiv \mathbf{0}$.
Consider $y\equiv \mathbf{0}$, according to the coupling relationship between $y_i$ and $(\xi_{1}^{1}, \cdots, \xi_{r_{1}}^{1}, \cdots, \xi_{1}^{m}, \cdots, \xi_{r_{m}}^{m})= \xi^\top$, we have $\xi \equiv \mathbf{0}$. Thus we only need to prove $z\equiv \mathbf{0}$, where $z$ satisfies
\begin{equation}
\begin{split}
\dot{z}&=N z+B_{a}'U_az,\\
\mathbf{0}&=Hz.
\end{split}
\end{equation}
Define $Hz=\hat{y}$ and $n_0-s=n$. The condition for completely observable becomes: if $\hat{y} \equiv \mathbf{0}$, there is $z \equiv \mathbf{0}$, which is also the condition for the following system to be observable:
\begin{equation}\label{eq11}
\begin{split}
\dot{z}&=N z+B_{a}'U_az,\\
\hat y&=Hz,
\end{split}
\end{equation}
where $z \in \mathbb{R}^{n}$ and $\hat{y} \in \mathbb{R}^{m}$ are state and output of the observability equivalent system respectively; $N \in \mathbb{R}^{n\times n}$ is system matrix and $B_a' \in \mathbb{R}^{n\times k}$ is input matrix. Now the problem becomes that the attacker wants to damage system observability using state feedback control $U_a$, while the defender wants to protect system observability by controlling $H$, which is illustrated in Figure \ref{fig:sub2}. 

\subsection{Game formulation}
Define an infinitely repeated game as a tuple $(N,A,J)$. $N=\{a,\,d\}$ is the set of players, where $a,\,d$ represent the attacker and the defender respectively. $A=A^{a}\oplus A^{d}$ is the action set of players, where $A^{a}=\{U_a\in \mathbb{R}^{k\times n}\}, A^{d}=\{H \in \mathbb{R}^{m\times n}\}$. $J=\{J^a,J^d\}$ is the utility function set of the attacker and the defender. Define the utility function of the attacker in epoch $i$ as
\begin{align*}
J^a(U_{a,i},H_i)= \operatorname{dim\,Ker}\!\! \left[\begin{array}{c}
H_i\\
H_i\left(N+B_{a}'U_{a,i}\right) \\
\vdots \\
H_i\left(N+B_{a}'U_{a,i}\right)^{n-1}
\end{array}\right]\overset{\triangle}{=}\operatorname{dim\,Ker}\Omega(U_{a,i},H_i),
\end{align*}
where $\Omega$ is the observability matrix of system \eqref{eq11}, $U_{a,i}$ and $H_i$ are actions of players in epoch $i$. 
As for the defender, 
$J^d(U_{a,i},H_i)=-J^a(U_{a,i},H_i)$.
Define the value function $\Phi$ in epoch $i$ as the dimension of unobservable subspace, i.e., 
\begin{equation}
\Phi(U_{a,i},H_i)=\operatorname{dim\,Ker}\Omega(U_{a,i},H_i).
\end{equation} 
Therefore, the defender aims to minimize the value function, while the attacker seeks to maximize it. Two players update actions asynchronously and an epoch is defined once a player acts. This asynchronous decision making between the defender and the attacker reflects industrial-control practice: the defender first configures the observations; only once data start flowing can the attacker tamper with them, naturally producing a defense-then-attack cycle. Repetitive games can effectively reflect the interaction between the defender and attacker, guiding both sides in designing strategies to gain an advantage in this tug-of-war relationship. 
\begin{remark}
The usual discounted-sum criterion for repeated games is not adopted here, because the stage utility is the discrete dimension of the unobservable subspace, which is a quantity that is inherently non-additive. Instead, the analysis focuses on the unobservable dimension that persists once the play converges (or settles into a limit cycle), rather than on any discounted accumulation of these dimensions.
\end{remark}
\section{The best response sets}
In this section, we will give derivations and algorithms to get the best response sets of the attacker and defender respectively. We assume that matrices $N, B_a'$ and players’s actions $U_a,H$ are public knowledge. 
\subsection{The best response set of the attacker}
The attacker aims to maximize the value function by controlling $U_a$, i.e., 
\begin{equation} \label{eq27}
U_a^{*}=\!arg \!\max _{U_a\in \mathbb{R}^{k\times n}}\! \operatorname{dim\,Ker}
\Omega(U_{a},H).
\end{equation}
Denote $\mathcal{V} \triangleq\operatorname{Ker}\, \Omega $. Then $\mathcal{V}$ is a controlled invariant subspace contained in $\operatorname{Ker}H$. Among all controlled invariant subspaces contained in $\operatorname{Ker}H$,
there is a maximal one denoted as $\mathcal{V}^{*}$, thus 
\begin{equation}\label{eq_BR1a}
\max _{U_a\in \mathbb{R}^{k\times n}} \operatorname{dim\,Ker}\Omega(U_{a},H)=\operatorname{dim}\mathcal{V}^*(H).
\end{equation}
Here follows a lemma to find $\mathcal{V}^{*}$.
\begin{lemma}  \cite{basile1969} \label{leV*}
Let \ $\mathcal{V}_{0}=\operatorname{Ker} H$\ and define, for $\ i=   0,1,2,$ \ldots, 
 \begin{equation} 
 \\\mathcal{V}_{i+1}=\left\{x \in \operatorname{Ker} H \mid N x \in \mathcal{V}_{i}+\operatorname{Im} B_a'\right\}. 
 \end{equation} 
Then $ \mathcal{V}_{i+1} \subset \mathcal{V}_{i}$. There exists $q \in \mathbb{R},\, q \leq \operatorname{dim\,} \mathcal{V}_{0}$, $\mathcal{V}_{q+1}=\mathcal{V}_{q}=\mathcal{V}^{*}$.
\end{lemma}
Classical results on maximal controlled invariant subspaces show that $U_a\in \mathbb{R}^{k\times n}$ maximizes the dimension of the unobservable subspace \emph{iff} it is a friend matrix of $\mathcal{V}^*$; formally,
\begin{equation*}
U_a^*\in \mathcal{F}(\mathcal{V}^*(H)),
\end{equation*}
where $\mathcal{F}(\mathcal{V}^*)$ is defined in Definition \ref{defi_1}. Algorithm 1 constructs such a friend matrix via the pseudo inverse. In most cases, the friend matrix is not unique; hence the pseudo inverse
provides a convenient explicit realization of \(U_a^{*}\).
The complete computational steps are summarized in Algorithm 1. 

\begin{algorithm}\label{alg:2}
    \caption{Attacker: maximization of the unobservable subspace}
    \begin{algorithmic}
    \STATE {\textbf{Input}  system matrices $N$, $B_a'$; the defender action $H$}\\
        \STATE {Set $\operatorname{Im}V_0=\operatorname{Ker}H.$}
        \\{Calculate $\operatorname{Im}Z_1=\left(\operatorname{Ker}\left[\begin{array}{c}V_0^{\top} \\B_a'^{\top}\end{array}\right]\right)^{\top}$, $\operatorname{Im}V_{1}=\operatorname{Ker}\left[\begin{array}{c}H \\Z_{1}N  \end{array}\right]$}.
        \STATE {Set $i=1$}
        \WHILE{$\operatorname{Im}V_i \neq \operatorname{Im}V_{i-1}$}
        \STATE Set $i=i+1$. \\{Calculate $\operatorname{Im}Z_{i}=\left(\operatorname{Ker}\left[\begin{array}{c}V_{i-1}^{\top} \\B_a'^{\top}\end{array}\right]\right)^{\top}$, $\operatorname{Im}V_{i}=\operatorname{Ker}\left[\begin{array}{c}H \\Z_{i}N\end{array}\right]$}.
        \ENDWHILE
        \\{Set $V=V_{i}$}.\\
        {Calculate $[X;U]=\text{pinv}\left(V\,\,\, B_a'\right) NV$, choose the last $(n-r)$ columns as $U$}.\\
        {Calculate $U_a^*=-U\text{pinv}(V)$}.\\
        {\textbf{return}  $U_a^*$}
    \end{algorithmic} 
\end{algorithm}  

\subsection{The best response set of the defender}
The defender aims to minimize the value function by controlling $H$, i.e.,
\begin{equation}\label{eq16}
{H^{*}}=\arg \min _{{H}\in \mathbb{R}^{m\times n}} \operatorname{dim\, Ker}\Omega(U_{a},H),
\end{equation}
which is equal to $\arg \max _{H\in \mathbb{R}^{m\times n}} \operatorname{dim\,Im}\Omega^T$,where 
\begin{equation*}
\Omega^T\!= [H^\top\!, \!\left(N+B_{a}'U_a\right)^\top \!H^\top,\!\cdots,\left(N+B_{a}'U_a\right)^{(n-1)\top}\!H^\top]
\end{equation*}
can be viewed as the controllability matrix of the system \eqref{eq11}'s dual system:
$\dot{\bar{z}}=\bar{A} \bar{z}+\bar{B} \bar{u}$,  
where $\bar{A}=(N+B_{a}'U_a)^\top, \bar{B}=H^\top$. Thus the best response set \eqref{eq16} becomes
$\bar{B}^*=\arg \max _{\bar{B}} \operatorname{dim\,Im}
\left[\begin{array}{lll}
\bar{B},\!\! &\bar{A} \bar{B},\cdots, \bar{A}^{n-1} \bar{B}\end{array}\!\!\right]$. The problem becomes how to choose $\bar{B}$ to make the dual system controllable. Find similar transformation matrix $T\in \mathbb{R}^{n\times n}$ which makes $\bar{A}$ become Jordan normal form $J$, i.e., $J=T^{-1}\bar{A}T$. And $\bar{B}$ becomes $\hat{B}=T^{-1} \bar{B}$.
For the above $J$, let its $l$ eigenvalues be:  $\lambda_{1}$ (algebraic multiplicity: $\sigma_{1}$, geometric multiplicity: $\alpha_{1}$), $\lambda_{2}$ (algebraic multiplicity: $\sigma_{2}$, geometric multiplicity: $\alpha_{2}$), $\cdots, \lambda_{l}$ (algebraic multiplicity: $\sigma_{l}$, geometric multiplicity: $\alpha_{l}$). Assume $\lambda_{i} \neq \lambda_{j}, \forall i \neq j$ and $\sigma_{1}+\sigma_{2}+\cdots+\sigma_{l}=n$. 
Thus we have
\begin{equation*} 
J=J\left(\lambda_{1}\right) \oplus J\left(\lambda_{2}\right) \oplus \cdots \oplus J\left(\lambda_{l}\right),\,
\hat{B}=\left[\hat{B}_{1}^{\top}, \hat{B}_{2}^{\top}, \cdots, \hat{B}_{l}^{\top}\right]^{\top},
\end{equation*}
where $\oplus$ is the direct sum of matrices,
\begin{equation*}
J(\lambda_{i})=
\begin{bmatrix}
 {J}_{1}(\lambda_{i}) & & & \\
&  {J}_{2}(\lambda_{i}) & & \\
& & \ddots & \\
& & &  {J}_{\alpha_{i}}(\lambda_{i})
\end{bmatrix},\, \hat{B}_{i}=\begin{bmatrix}
 \hat{ {B}}_{i 1} \\
\hat{ {B}}_{i 2} \\
\vdots \\
\hat{ {B}}_{i \alpha_{i}}   
\end{bmatrix},  
\end{equation*}
for $i=1,2,\cdots,l$, where
\begin{equation*}
J_{k}(\lambda_{i})=\begin{bmatrix}
\lambda_{i} & 1 & & & \\
& \lambda_{i} & 1 & & \\
& & \ddots & \ddots & \\
& & & \ddots & 1 \\
& & & & \lambda_{i}
\end{bmatrix},\,\hat{ {B}}_{i k}=\begin{bmatrix}
\hat{ {b}}^{ i k}_{1} \\
\hat{ {b}}^{ i k}_{2} \\
\vdots  \\
\hat{ {b}}^{ i k}_{r_{ik}-1}\\
\hat{ {b}}^{ i k}_{r_{ik}}
\end{bmatrix},
\end{equation*}  
for $k=1, 2,\cdots, \alpha_{i}$, where $J_{k}(\lambda_{i}) \in \mathbb{R}^{r_{i k} \times r_{i k}}, \hat{ {B}}_{i k} \in \mathbb{R}^{r_{i k} \times m}$, $\sum_{k=1}^{\alpha_{i}} r_{i k}=\sigma_{i}$.
Recall the following Lemma which gives the controllability condition in the Jordan‐form representation.
\begin{lemma}\cite{chen1984linear}\label{jordan} 
For the Jordan normal form of linear system, the necessary and sufficient condition for complete controllability is
\begin{equation} 
\operatorname{rank}\left[
\hat{b}_{r_{i1}}^{i1\top}, \hat{b}_{r_{i2}}^{i2\top}, 
\cdots,
\hat{b}_{r_{i\alpha_i}}^{i\alpha_i\top}\right]=\alpha_{i}, \forall i=1,2, \cdots,l,
\end{equation} 
where $\alpha_{i}$ is the geometric multiplicity of eigenvalue $\lambda_{i}$ of matrix $\bar{A}$. This means the last rows of $\hat{ {B}}_{i 1} ,\hat{ {B}}_{i 2} ,\cdots,
\hat{ {B}}_{i \alpha_{i}}$ are linearly independent.
\end{lemma} 

However, if there is a geometric multiplicity $\alpha_{i}$ exceeding the number of columns of $\hat{B}$, i.e., $\exists \alpha_i>m$, it is impossible for the system to be fully controllable. In this case, we give a construction of $\hat{B}$ which maximizes the controllable subspace. 
\begin{proposition}\label{pro4}
For $\hat{B}$, we assume that the last rows of $\hat{ {B}}_{i 1} ,\hat{ {B}}_{i 2} ,\cdots,
\hat{ {B}}_{i \alpha_{i}}$ satisfy
\begin{equation} \label{eq20}
\left[\!\begin{array}{c}
\hat{b}_{r_{i1}}^{i1} \\
\hat{b}_{r_{i2}}^{i2} \\
\vdots \\
\hat{b}_{r_{i\alpha_i}}^{i\alpha_i}
\end{array}\!\right]\!\!=\!\!
\left\{ 
\begin{array}{cc}
\left[I_{\alpha_{i}}\ \textbf{0}_{\alpha_{i} \times\left(m-\alpha_{i}\right)}\right],\, \text{for} \,\{i \mid m \geq \alpha_{i}\},\\
\\
\left[\delta_{\alpha_i}^{j_1},\delta_{\alpha_i}^{j_2},\cdots,\delta_{\alpha_i}^{j_m}\right], \,\text{for} \,\{i \mid m < \alpha_{i}\},
\end{array}
\right.
\end{equation}
where \,$\delta_{\alpha_i}^{j_p}=\operatorname{Col}_{j_p}\left(I_{\alpha_i}\right)$, $\{j_1,j_2,\cdots,j_m\} \subset \{1,2,\cdots,\alpha_i\}$ are the subscripts of $m$-th largest Jordan blocks for eigenvalue $\lambda_i$ and other rows of $\hat{ {B}}_{i 1} ,\hat{ {B}}_{i 2} ,\cdots ,
\hat{ {B}}_{i \alpha_{i}}$ are zero rows. Then the controllable subspace is maximized.

\end{proposition}
\begin{proof}
According to Lemma \ref{jordan}, we can easily derive that $\hat{B}\in \mathbb{R}^{n\times m}$ maximizes the dimension of the controllable subspace if and only if the last rows of $\hat{ {B}}_{i 1} ,\hat{ {B}}_{i 2} ,\cdots,
\hat{ {B}}_{i \alpha_{i}}$ satisfy
\begin{equation}\label{eq19}
\left\{ 
\begin{array}{cc}
\operatorname{rank}\left[
\hat{b}_{r_{i1}}^{i1\top}, \hat{b}_{r_{i2}}^{i2\top}, 
\cdots,
\hat{b}_{r_{i\alpha_i}}^{i\alpha_i\top}\right]=\alpha_{i} ,\,\text{for} \,\{i \mid m \geq \alpha_{i}\},\\
\\
\operatorname{rank}\left[
\hat{b}_{r_{ij_1}}^{ij_1\top}, \hat{b}_{r_{ij_2}}^{ij_2\top}, 
\cdots,
\hat{b}_{r_{ij_m}}^{ij_m\top}\right]=m ,\,  \text{for} \,\{i \mid m < \alpha_{i}\},
\end{array}
\right.
\end{equation}
where $\{j_1,j_2,\cdots,j_m\} \subset \{1,2,\cdots,\alpha_i\}$ are the subscripts of $m$-th largest Jordan blocks for eigenvalue $\lambda_i$.
For $m \geq \alpha_{i}$, $\operatorname{rank}\left[I_{\alpha_{i}}\ \textbf{0}_{\alpha_{i} \times\left(m-\alpha_{i}\right)}\right]=\alpha_i$; for $m < \alpha_{i}$, $\operatorname{rank}\left[\delta_{\alpha_i}^{j_1},\delta_{\alpha_i}^{j_2},\cdots,\delta_{\alpha_i}^{j_m}\right]=m$. Thus the construction of $\hat{B}$ in equation \eqref{eq20} maximizes the dimension of controllable subspace. This completes the proof.
\end{proof}
With $\hat{B}$ which satisfies Proposition \ref{pro4}, the best response for the defender is
\begin{equation*}
H^* \in \{\hat{B}^{\top} T^{\top}\|\hat{B} \text{ satisfies } \eqref{eq19} \},
\end{equation*}
where $T$ is the transformation matrix making $(N+B_{a}'U_a)^\top$ become Jordan normal form. The value function with ${H}^*\in \mathbb{R}^{m\times n}$ is 
\begin{equation}\label{eq_value_H}
\left\{ 
\begin{aligned}
&\min _{{H}\in \mathbb{R}^{m\times n}} \operatorname{dim\, Ker}\Omega(U_{a},H)=\operatorname{dim\, Ker}\Omega(U_{a},H^*)\!=\!n\!-\!\underset{{q}\in \mathcal{I}}{\min} (\sum_{p=1}^{m}r_{qj_p }),\,\text{for } \mathcal{I} \neq \emptyset,\\
&\min _{{H}\in \mathbb{R}^{m\times n}} \operatorname{dim\, Ker}\Omega(U_{a},H)=\operatorname{dim\, Ker}\Omega(U_{a},H^*)=0,\,  \text{for } \mathcal{I} = \emptyset,
\end{aligned}
\right.
\end{equation}
where $\mathcal{I}=\{i \mid \alpha_{i}>m\}$ includes subscripts for the eigenvalues of $(N+B_a'U_a)$ whose geometric multiplicity is larger than $m$, $\{r_{qj_1},r_{qj_2},\cdots,r_{qj_m}\}$ are the dimensions of $m$-th largest Jordan blocks for eigenvalue $\lambda_q$.
The calculation steps of $H^*$ based on Proposition \ref{pro4} is summarized in Algorithm 2.
\begin{algorithm}[H] \label{alg:2}
    \caption{Defender: minimization of the unobservable subspace}
\begin{algorithmic}
\STATE {\textbf{Input}  system matrices $N$, $B_a'$; the attacker action $U_a$}\\
\STATE {Set $\bar{A}=(N+B_{a}'U_a)^\top.$}
  \\{Compute the Jordan normal form of $\bar{A}:J=T^{-1}\bar{A}T$, whose distinct eigenvalues are $\lambda_1,\dots,\lambda_l$}.
   \FOR{$i = 1$ \TO $l$}
    \STATE{$\alpha_i\gets$ geometric multiplicity of $\lambda_i$; }
    \STATE{$\{r_{i1},r_{i2},\cdots,r_{i\alpha_i}\}\gets$ dimensions  of Jordan blocks of $\lambda_i$}; 
    \STATE{$\{j_1,j_2,\cdots,j_m\}\gets$ the subscripts of $m$-th largest Jordan blocks of $\lambda_i$.}
    \ENDFOR
        \\{Compute $\hat B$\ according to Corollary \ref{pro4}.}
         \\{$H^*=\hat{B}^{\top} T^{\top}$}.\\
        {\textbf{return}  $H^*$}
    \end{algorithmic}
\end{algorithm}

\section{Equilibrium analysis}
Based on the above best response sets, we next give the equilibrium of the game considering one-step and two-step optimality respectively in subsections 4.1 and 4.2. Then in subsections 4.3 and 4.4, game outcomes and equilibrium characterization are refined, which are suitable for both one-step or two-step optimality.
Finally, three key insights of the observability-adversarial game are summarized in subsection 4.5.
\subsection{Game based on one-step optimality}
In the one‐step optimality formulation of the repeated game, each player maximizes only the immediate utility at each stage. We denote the resulting best response sets for the attacker and the defender by $BR1^{a}$ and $BR1^{d}$, i.e.,
\begin{equation}
BR1^{a}(H)= \arg \max _{{U_a}\in \mathbb{R}^{k\times n}} \operatorname{dim\,Ker}\Omega(U_a,H),\,
 BR1^{d}(U_a)= \arg \min _{{H}\in \mathbb{R}^{m\times n}} \operatorname{dim\,Ker}\Omega(U_a,H).
\end{equation}
The solutions of these best response sets have been discussed in subsections 3.1 and 3.2, i.e.,  
\begin{equation}
BR1^{a}(H) =\mathcal{F}(\mathcal{V}^*(H)),
\end{equation}
where $\mathcal{F}(\mathcal{V}^*)$ is the friend matrix set of $\mathcal{V}^*$, $\mathcal{V}^*$ is the maximal controlled invariant subspace in $\operatorname{Ker} H$; 
\begin{equation}
BR1^{d}(U_a)=\{\hat{B}^{\top} T^{\top}\|\hat{B} \text{ satisfies } \eqref{eq19} \},
\end{equation}
where $T$ is the transformation matrix making $(N+B_{a}'U_a)^\top$ become Jordan normal form. Figure 2 shows the sequence of actions. 
\begin{center}
  \centerline
  {\includegraphics[scale=0.5]{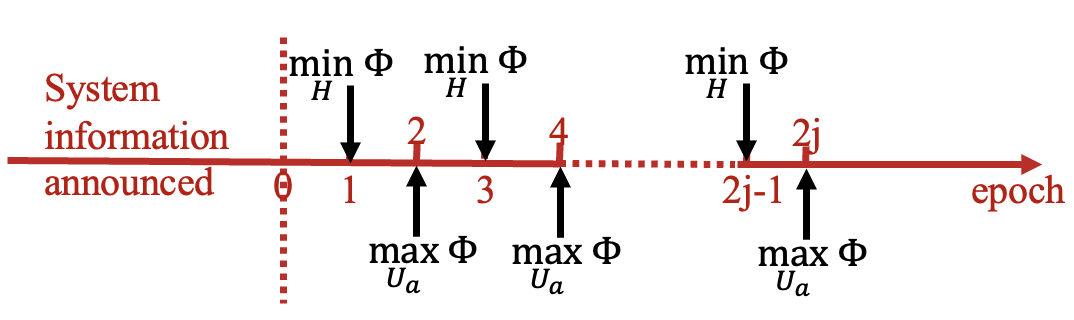}}
\centering{\small {\bf Figure 2}\ \ Sequence of actions considering one-step optimality. }
\end{center}

Both best response sets are not single-valued maps and Algorithm 1 (or Algorithm 2) only chooses a special $U_a$ (or $H$) in $BR1^a$ (or $BR1^d$). If a player chooses a different action in the best response set, the game result will be different. 

Define the Nash equilibrium based on one-step optimality as follows.
\begin{definition} \label{def3}
The strategy profile $(U_a^{*}, H^{*})$ is said to be the Nash equilibrium (NE) of the one-step optimal game, if $U_a^{*} \in BR1^{a}(H^{*}), H^{*} \in BR1^{d}(U_a^{*})$.
\end{definition}
Then we give a necessary and sufficient condition for the Nash equilibrium.
\begin{theorem}[One-step optimality NE criterion] \label{the1}
The strategy profile $(U_a^{*},H^{*})$ is a Nash equilibrium of the one-step optimal game if and only if 
\begin{equation}\label{equ_condition}
\min _{{H}\in \mathbb{R}^{m\times n}}\!\operatorname{dim\, Ker}\Omega \,(U_a^{*},\!H)=\operatorname{dim}\mathcal{V}^*(H^{*}).
\end{equation}
\end{theorem}
\begin{proof}
(Sufficiency) We prove $U_a^{*} \in BR1^{a}(H^{*}), H^{*} \in BR1^{d}(U_a^{*})$ by establishing its contrapositive. 
If $U_a^{*} \notin BR1^{a}(H^{*})$, there is $\operatorname{dim\, Ker}\Omega(U_a^{*},H^{*})<\operatorname{dim}\mathcal{V}^*(H^{*})$. Because there is $\min _{{H}\in \mathbb{R}^{m\times n}}\!\operatorname{dim\, Ker}\Omega \,(U_a^{*},\!H) \leq \operatorname{dim\, Ker}\Omega(U_a^{*},H^{*})$, this contradicts \eqref{equ_condition}. 
If $H^{*} \notin BR1^{d}(U_a^{*})$, there is $\min _{{H}\in \mathbb{R}^{m\times n}}\!\operatorname{dim\, Ker}\Omega \,(U_a^{*},\!H) < \operatorname{dim\, Ker}\Omega(U_a^{*},H^{*})$. Because there is $\operatorname{dim\, Ker}\Omega(U_a^{*},H^{*})$ $\leq \operatorname{dim}\mathcal{V}^*(H^{*})$, this contradicts \eqref{equ_condition}. Thus there must be $U_a^{*} \in BR1^{a}(H^{*}), H^{*} \in BR1^{d}(U_a^{*})$ and $(U_a^{*},H^{*})$ is a Nash equilibrium.

(Necessity) Because $U_a^{*} \in BR1^{a}(H^{*})$, there is 
\begin{equation}\label{eq_nece1}
\operatorname{dim\, Ker}\Omega(U_a^{*},H^{*})\,=\operatorname{dim}\mathcal{V}^*(H^{*}).
\end{equation}
Since $H^{*} \in BR1^{d}(U_a^{*})$, we have 
\begin{equation}\label{eq_nece2}
\min _{{H}\in \mathbb{R}^{m\times n}}\!\operatorname{dim\, Ker}\Omega \,(U_a^{*},\!H) = \operatorname{dim\, Ker}\Omega(U_a^{*},H^{*}).
\end{equation}
Combining \eqref{eq_nece1} and \eqref{eq_nece2}, we can get \eqref{equ_condition}, which completes the proof. 
\end{proof}
Theorem \ref{the1} is effective in determining when the Nash equilibrium is reached in dynamic game process, as it only requires verifying whether the value function $\operatorname{dim\, Ker}\Omega \,(U_a^{*},\!H)$ after the defender chooses $H^{*}$ is equal to $\operatorname{dim}\mathcal{V}^*(H^{*})$.

\subsection{Game based on two-step optimality}
Since the one-step best response sets, are multi-valued, different preferences for choices in \(BR1^{a}\) and \(BR1^{d}\) can steer the play toward different result.  
For any \(H\in BR1^{d}\), the induced value \(\dim\mathcal{V}^{*}\) varies, thereby tightening or relaxing the upper bound of \(\max_{U_a}\dim\operatorname{Ker}\Omega(U_{a},H)\).  
Conversely, every \(U_{a}\in BR1^{a}\) alters the maximum geometric multiplicity among all eigenvalues of \(N+B_{a}^{\top}U_{a}\), which fixes the lower bound of \(\min_{H}\dim\operatorname{Ker}\Omega(U_{a},H)\).

To capture these continuation effects, we impose a two-step optimality criterion: within \(BR1^{a}\) or \(BR1^{d}\), each player selects the action that maximizes the loss of the opponent in the subsequent stage.
The attacker’s two-step best response set is  
\begin{equation}\label{BR2_a}
BR2^{a}(H)=\arg \max _{U_{a} \in BR1^{a}(H)}\min _{H\in \mathbb{R}^{m\times n}} \operatorname{dim\, Ker} \Omega\,(U_a,H) =\arg \max _{U_{a} \in BR1^{a}(H)} [\text{MGM}(N+B_a'U_a)].
\end{equation}
where $\text{MGM}(.)$ is the maximum geometric multiplicity among all eigenvalues of the argument matrix, that is, the largest dimension of any eigen-space. This result can be derived from \eqref{eq_value_H}. 
The two-step best response set of the defender is  
\begin{equation}\label{BR2_d}
BR2^{d}(U_{a})=arg \min _{H \in BR1^{d}(U_{a})}\max _{U_{a}\in \mathbb{R}^{k\times n}} \operatorname{dim\, Ker} \Omega\,(U_a,H)=arg \min _{H \in BR1^{d}(U_{a})} \operatorname{dim\, }\mathcal{V}^*(H),
\end{equation}
which can be derived from \eqref{eq_BR1a}. 

For comparison, define $BR2X^{a}\triangleq \arg \max _{U_{a} \in \mathbb{R}^{k\times n}} [\text{MGM}(N+B_a'U_a)]$, $BR2X^{d}\triangleq \arg \,$ $ \min _{H \in \mathbb{R}^{m\times n}} \operatorname{dim\, }\mathcal{V}^*(H)$.
Neither \(BR2X^{a}\) nor \(BR2X^{d}\) is a valid two‐step best response, since each ignores the requirement to optimize the current period’s value function. Instead, they
correspond exactly to the Stackelberg solutions under two different leadership orders: when the defender leads, it commits to \(BR2X^{d}\) and the attacker
responds with \(BR1^{a}\); when the attacker leads, it chooses
\(BR2X^{a}\) and the defender replies with \(BR1^{d}\).

Figure 3 shows the sequence of actions following two-step optimality. 
\begin{center}
  \centerline
  {\includegraphics[scale=0.4]{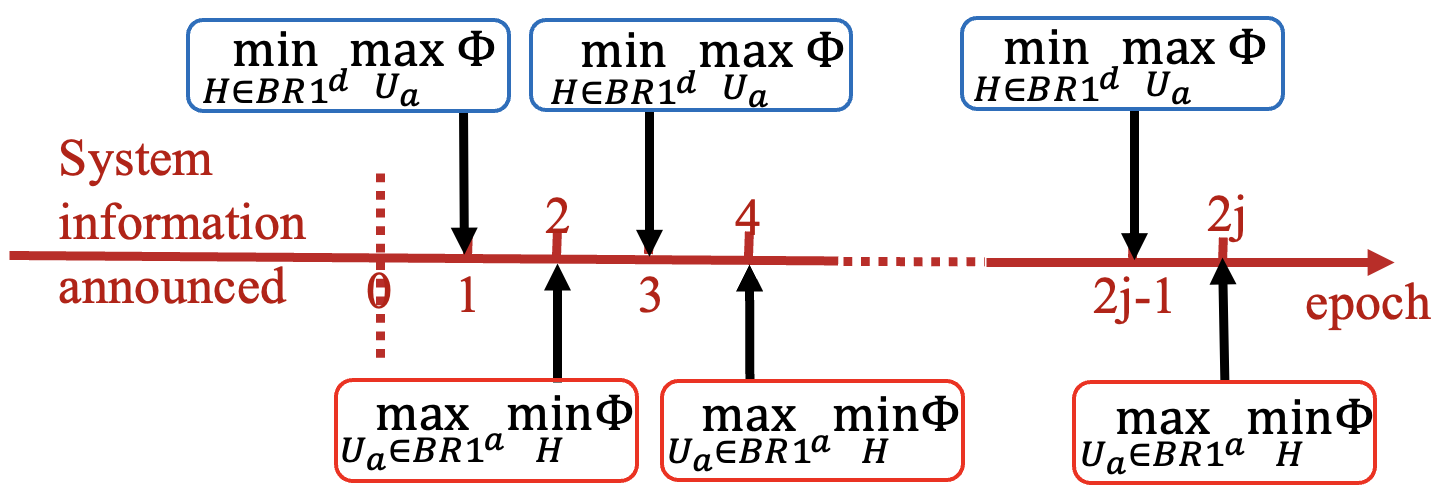}}
\centering{\small {\bf Figure 3}\ \ Sequence of actions following two-step optimality. }
\end{center}

Then we define the Nash equilibrium when two players consider two-step optimality.
\begin{definition}
The strategy profile $(U_a^{*}, H^{*})$ is said to be a Nash equilibrium of the two-step optimal game if $ U_a^{*} \in BR2^{a}(H^{*}), H^{*} \in BR2^{d}(U_a^{*})$.
\end{definition}

In order to get the condition for Nash equilibrium under two-step optimality, we give the following lemma.
\begin{lemma}\label{le4}
For the same $N \in \mathbb{R}^{n\times n}, B_a'\in \mathbb{R}^{n\times k}$, $\forall U_a\in \mathbb{R}^{k\times n}$, there is 
\begin{equation}
\min _{H \in \mathbb{R}^{m\times n}} \operatorname{dim\, }\mathcal{V}^*(H)\geq \min _{H \in \mathbb{R}^{m\times n}} \operatorname{dim\, Ker}\Omega(U_a,H).
\end{equation}
\end{lemma}
\begin{proof}
$\forall \, H_0\in \mathbb{R}^{m\times n}$, the maximal controlled invariant subspace in $\operatorname{Ker} H_0$ is larger than $(N+B_a'U_a)$-invariant subspace in $\operatorname{Ker} H_0$, i.e., $\operatorname{dim\,}\mathcal{V}^*(H_0) \geq \operatorname{dim\, Ker}\Omega(U_a,H_0).$ Thus $\forall H_1 \in arg \min _{H \in \mathbb{R}^{m\times n}} \operatorname{dim\, }\mathcal{V}^*(H)$, $\operatorname{dim\,}\mathcal{V}^*(H_1) \geq \operatorname{dim\, Ker}\Omega(U_a,H_1)$. And $\forall H_2 \in \arg \min _{{H}\in \mathbb{R}^{m\times n}}\,$ $ \operatorname{dim\,Ker}$ $\Omega(U_a,H)$, there is $\operatorname{dim\, Ker}\Omega(U_a,H_1)\geq \operatorname{dim\, Ker}\Omega(U_a,H_2)$. Thus $\operatorname{dim\,}\mathcal{V}^*(H_1)\!\geq\! \operatorname{dim\, Ker}\Omega(U_a,H_2)$, which completes the proof.
\end{proof}
Building on this lemma, we demonstrate that the Nash-equilibrium criterion stated in Theorem \ref{the1} also holds when both players adopt a two-step optimality framework. Although the equilibrium condition is identical for the one-step and two-step settings, the underlying mathematical reasoning are distinct.
\begin{theorem}[Two-step optimality NE criterion] \label{the2}
The strategy profile $(U_a^{*},H^{*})$ is a Nash equilibrium when two players consider two-step optimality if and only if, 
\begin{equation}
\min _{H\in \mathbb{R}^{m\times n}}\operatorname{dim\, Ker}\Omega\,(U_a^{*},H)=\operatorname{dim\,}\mathcal{V}^*(H^{*}).
\end{equation}
\end{theorem}
\begin{proof} 
(Sufficiency) According to Theorem \ref{the1}, $H^{*} \in BR1^d(U_a^{*})$, $U_a^{*} \in BR1^a(H^{*})$ and
\begin{equation}\label{eq_the2}
\min _{H\in \mathbb{R}^{m\times n}}\operatorname{dim\, Ker}\Omega\,(U_a^{*},H)=\operatorname{dim\,}\mathcal{V}^*(H^{*})=\operatorname{dim\, Ker}\Omega\,(U_a^{*},H^*)\overset{\triangle}{=}\gamma.
\end{equation}
We need to further prove $H^{*} \in BR2^d(U_d^{*})$ and $U_a^{*} \in BR2^a(H^{*})$. Firstly, according to Lemma \ref{le4} and \eqref{eq_the2}, $\forall H' \in BR1^d(U_a^{*})$ we have
\begin{equation}
\operatorname{dim\,}\mathcal{V}^*(H') \geq \min _{H \in \mathbb{R}^{m\times n}} \operatorname{dim\, }\mathcal{V}^*(H)\geq \min _{H \in \mathbb{R}^{m\times n}} \operatorname{dim\, Ker}\Omega(U_a^*,H)=\gamma.
\end{equation}
By \eqref{eq_the2}, $\operatorname{dim\,}\mathcal{V}^*(H^{*})=\gamma$, which reaches the lower bound of the $\operatorname{dim\,}\mathcal{V}^*(H\in BR1^{d}(U_a^{*}))$. Thus $H^* \in arg \min _{H \in BR1^{d}(U_{a}^{*})} \operatorname{dim\, }\mathcal{V}^*(H)=BR2^{d}(U_a^{*})$. Secondly, by Lemma \ref{le4} and \eqref{eq_the2}, $\forall U_a' \in BR1^a(H^{*})$, there is
\begin{equation}
\min _{H \in \mathbb{R}^{m\times n}} \operatorname{dim\, Ker}\Omega(U_a',H) \leq \min _{H \in \mathbb{R}^{m\times n}} \operatorname{dim\, }\mathcal{V}^*(H) \leq \operatorname{dim\,}\mathcal{V}^*(H^*) =\gamma.
\end{equation}
By \eqref{eq_the2}, $\min _{H \in \mathbb{R}^{m\times n}} \operatorname{dim\, Ker}\Omega(U_a^*,H)=\gamma$, which reaches the upper bound of $\min _{H \in \mathbb{R}^{m\times n}} \operatorname{dim\, Ker}$ $\Omega(U_a',H)$. Thus $U_a^* \in \arg \max _{U_{a} \in BR1^{a}(H)}\min _{H\in \mathbb{R}^{m\times n}} \operatorname{dim\, Ker} \Omega\,(U_a,H)=BR2^{a}(H^*)$. Thus $(U_a^{*},H^{*})$ is the Nash equilibrium of the two-step optimal game.
\\(Necessity) Since $BR2^d \subseteq BR1^d$ and $BR2^a \subseteq BR1^a$, we have $H^{*} \in BR1^d(U_a^{*})$, $U_a^{*} \in BR1^a(H^{*})$. The rest proof is the same as Theorem \ref{the1}.
\end{proof}

Theorem \ref{the2} reveals that, although the two-step best response sets are a subset of the one-step best response sets and the strategy sets of both players differ, the resulting equilibrium points coincide.
\begin{remark}
In this paper, the analysis is limited to one-step and two-step optimality. Extending the horizon to three or more steps introduces no additional strategic content, because the same player acts in both epochs~1 and~3, making the strategic situation in epoch~3 identical to that in epoch~1. Hence, higher‐order optimality criteria can be omitted without loss of
generality, and the extended-horizon study in Section~4.2 (two-step optimality) exhausts all non-trivial multi-step cases.
\end{remark}

\subsection{Game outcome analysis}
This subsection gives analysis for game results in dynamic game process that are consistent with both one-step and two-step optimality game. First, we examine a degenerate scenario in which the game possesses an equilibrium where the defender holds an absolute advantage. The following lemma states sufficient conditions under which this situation occurs.
\begin{theorem}[Defender-dominated NE]\label{th3}
Let $N \in \mathbb{R}^{n\times n}$ and $B_a' \in \mathbb{R}^{n\times k}$. Assume there exists $H^{*}\in\mathbb{R}^{m\times n}$ such that either of the following two cases holds:
\\Case 1: $(n-m)\geq k. \text { If } \operatorname{Im}B_a' \subseteq \operatorname{Ker} H^*$ and there is no nontrivial $N$-invariant subspace contained in $\operatorname{Ker} H^*$; 
\\Case 2: $(n-m)<k$. If $\operatorname{Ker} H^* \subseteq \operatorname{Im}B_a'$ and any vector in  $\operatorname{Ker} H^*$ does not belong to $N$-invariant subspace contained in $\operatorname{Im}B_a'$.
Then 
\begin{equation}
\operatorname{dim}\mathcal{V}^*(H^{*})=0,
\end{equation}
and, for every $U_a \in \mathbb{R}^{k\times n}$, the strategy profile $(U_a,H^{*})$ is a Nash equilibrium.
\end{theorem}
\begin{proof}
$(1)\text { When } (n-m)\geq k. \text { Since } \operatorname{Im}B_a' \subseteq \operatorname{Ker} H^*, H^*B_a'=0. \text { Then } \forall v \neq 0 \in \operatorname{Ker} H^*$, $\forall U_a,  H^*(N+B_a'U_a)v=H^* N v \neq 0.$ Thus $(N+B_a'U_a)v \notin \operatorname{Ker} H^*$, which means controlled invariant subspace contained in $\operatorname{Ker} H^*$ is $\textbf{0}$. Thus dim $\mathcal{V}^*(H^{*})=0.$
(2)$\text { When }(n-m)<k. \forall v \neq 0 \in \operatorname{Ker} H^*, v \in \operatorname{Im} B_a'$, because $Nv \notin \operatorname{Im} B_a', (N+B_a'U_a)v\notin \operatorname{Im} B_a'$. Since $\operatorname{Ker} H^* \subseteq \operatorname{Im}B_a',(N+B_a'U_a)v \notin \operatorname{Ker} H^*$, which means controlled invariant subspace contained in $Ker H^*$ is $\textbf{0}$. Thus dim $\mathcal{V}^*(H^{*})=0$.

Furthermore, when dim $\mathcal{V}^*(H^{*})=0$, $\forall U_a \in \mathbb{R}^{k\times n}$, $\operatorname{dim\, Ker}\Omega$ $(U_a,H^{*})=\operatorname{dim\,}\mathcal{V}^*(H^{*})=0$. Thus $H^{*} \in \arg \min _{H\in \mathbb{R}^{m\times n}}\operatorname{dim\, Ker}\Omega$ $\,(U_a,H)$. Plus, there is $\min _{H\in \mathbb{R}^{m\times n}}\operatorname{dim\, Ker}\Omega\,(U_a,H)=\operatorname{dim\,}\mathcal{V}^*$ $(H^{*})=\operatorname{dim\, Ker}\Omega(U_a,H^{*})$. According to Theorem \ref{the1} and Theorem \ref{the2}, $(U_a,H^{*}), \forall U_a \in \mathbb{R}^{k\times n}$ is an equilibrium whenever two players consider one-step or two-step optimality, which completes the proof. 
\end{proof}

When the conditions in Theorem~\ref{th3} hold, the defender can choose a matrix \(H^{*}\) such that \(\mathcal{V}^*(H^{*})=\{0\}\).
In this defender-dominated situation the attacker can no longer influence
the system’s observability, and the game settles at a trivial equilibrium
completely controlled by the defender.
To exclude this degenerate case, the remainder of the paper focuses on the
non-trivial regime $\min _{H\in \mathbb{R}^{m\times n}}\operatorname{dim\,}\mathcal{V}^*>0$.
The following theorem provides a convenient sufficient condition under which this inequality is guaranteed.
\begin{corollary}[Non-degenerate condition]\label{pro10}
Consider system~\eqref{eq11}.  Assume that
\(B_a'\in\mathbb{R}^{n\times k}\) and \(H\in\mathbb{R}^{m\times n}\) are both full-column-rank.  If
\begin{equation}\label{eq_non_degene}
  \max \bigl\{k,\operatorname{MGM}(N_{\bar r})\bigr\}>m,
\end{equation}
where \(\operatorname{MGM}(N_{\bar r})\) denotes the maximum geometric
multiplicity among all eigenvalues of the uncontrollable block
\(N_{\bar r}\) of the state matrix \(N\), then
\begin{equation}
  \min_{H\in\mathbb{R}^{m\times n}}\dim\mathcal{V}^*(H)>0.
\end{equation}
\end{corollary}
\begin{proof}
Split the system into controllable and uncontrollable parts:
\(N=\text{diag}(N_r,N_{\bar r})\),
\(B_a'^{\!\top}=(B_{ar}'^{\!\top},0)\), and \(H=(H_r,H_{\bar r})\).
Hence \(\mathcal{V}^*=\mathcal{V}_r^*\oplus\mathcal{V}_{\bar r}^*\) with
\(\mathcal{V}_r^*,\mathcal{V}_{\bar r}^*\) determined by the subsystems
\((N_r,B_{ar}',H_r)\) and \((N_{\bar r},0,H_{\bar r})\), respectively. Consider the following two cases.
\textbf{Case 1}: If $k\geq \text{MGM}({N}_{\bar{r}})$, there is $k > m$ by \eqref{eq_non_degene}. For controllable part $(N_{r},B'_{ar},H_{r})$, the geometric multiplicity of $(N_{r}+B'_{ar}U_{ar})'s$ eigenvalue $\lambda$ is $\{n-\text{rank}[(N_{r}+B'_{ar}U_{ar})-\lambda I]\}.$ According to PBH controllability criterion, for the controllable part we have $\forall \lambda \in \sigma(N_{r}+B'_{ar}U_{ar})$, rank$[(N_{r}\!+\!B'_{ar}U_{ar}\!-\!\lambda I)\,\,\, B'_{ar}]=n$. Thus rank$[(N_{r}\!+\!B'_{ar}U_{ar}\!-\!\lambda I)]\geq (n-k)$ and $\text{MGM}(N_{r}\!+\!B'_{ar}U_{ar})=k.$ If $k>m$, i.e., $\text{MGM}(N_{r}\!+\!B'_{ar}U_{ar})>m$, according to Proposition \ref{pro4}, $\mathcal{I}=\{i \mid \alpha_{i}>m\} \neq \emptyset, \underset{{H_r}\in \mathbb{R}^{m\times n}}{\min} \operatorname{dim\,Ker}\Omega(U_{ar},H_r)=n-\underset{{q}\in \mathcal{I}}{\min} (\sum_{p=1}^{m}r_{qj_p })>0.$ According to Lemma \ref{le4}, $\min _{H_r \in \mathbb{R}^{m\times n}} \operatorname{dim\, }\mathcal{V}_r^*\geq \min _{H_r \in \mathbb{R}^{m\times n}} \operatorname{dim\, Ker}\Omega(U_{ar},H_r)>0$.
Thus $\min _{{H}_r} \operatorname{dim\, }$ $\mathcal{V}_r^*>0$. 
\textbf{Case 2}: If $\text{MGM}({N}_{\overline{r}})>k$, there is $\text{MGM}({N}_{\overline{r}})>m$ by \eqref{eq_non_degene}. Since $\mathcal{V}^*_{\bar{r}}=\operatorname{Ker}\Omega(N_{\bar{r}},0,H_{\bar{r}})$ for uncontrollable subspace, if $\text{MGM}({N}_{\overline{r}})>m$, extending Lemma \ref{jordan} to the dual system, we have $\min _{{H}_{\overline{\mathrm{r}}}} \operatorname{dim\, Ker}\Omega>0$. Thus $\min_{{H}_{\overline{\mathrm{r}}}}\operatorname{dim\, }\mathcal{V}_{\bar{r}}^*>0$. To conclude,
in either case, there is 
$\operatorname{dim\, }\mathcal{V}^*=\operatorname{dim\, }\mathcal{V}_r^*+\operatorname{dim\, }\mathcal{V}_{\bar{r}}^*>0,$ which completes the proof. 
\end{proof}
Next, under the non-degenerate condition, we analyze the outcome of the game. Since the best response sets for both players are not single-valued mappings, we impose the following assumptions on strategy selection of both players:

\textbf{Assumption 1:} The attacker (or the defender) prefers to keep the action unchanged if the last action also belongs to the best response set in this epoch.

\textbf{Assumption 2:} Without violating Assumption 1, if the best response set is the same in different game epochs, the attacker (or defender) consistently chooses the same action as the first time.

To further illustrate the above assumptions, Figure 4 represents an example of the evolution of the best response sets and actions for one player, who updates its action every two epochs. 
\begin{center}
  \centerline
  {\includegraphics[scale=0.4]{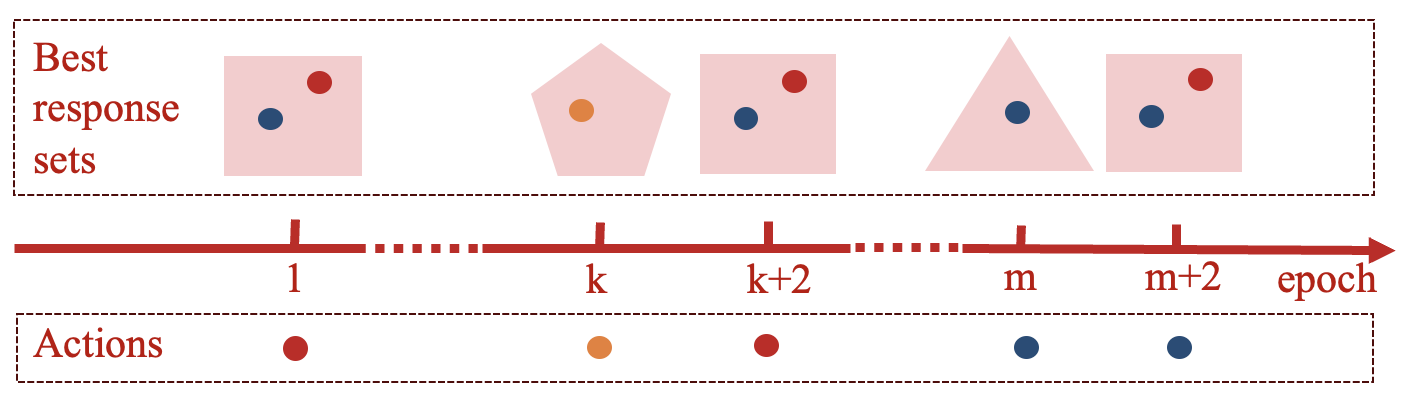}}
\centering{\small {\bf Figure 4}\ \ An illustration of how best response sets and chosen actions evolve over game epochs: shapes above the axis denote the best response set of each epoch, while colored points below indicate the action selected at that epoch. }
\end{center}

In epoch 1, the player chooses an action (the red point) in the best response set (the square pattern). Then in epoch $k+2$, the player has the same best response set as that in epoch 1 and chooses the same action as that in epoch 1 according to Assumption 2. Afterwards, in epoch $m+2$, though the player has the same best response set as that in epoch 1, it keeps the action unchanged (the blue point), because the last action (the blue point) also belongs to the best response set in this epoch and Assumption 1 is prioritized over Assumption 2. 

Assumptions 1 and 2 frequently manifest in real-world applications. Assumption 1 can save energy used in changing actions, while Assumption 2 fixes the rule to select a solution in a multi-valued map (such as choosing the one with the smallest norm length). 
Based on these assumptions, we can easily get the following theorem which represents two special game results.

\begin{theorem}[Game results analysis]
\label{cor7}
Under Assumptions~\textbf{1}–\textbf{2}, there are two possible outcomes in the infinitely repeated game when considering either one-step or two-step optimality:
\begin{enumerate}[label=\textnormal{(\roman*)}]
  \item Lock mode.  
        \(\exists\, l, \gamma \in\mathbb R\), \(\forall i\ge l\)  such that \(U_{a,i}=U_{a,l}\), \(H_i=H_l\) and
        \(\dim\operatorname{Ker } \Omega(U_{a,i},H_i)=dim\operatorname{Ker } \Omega(U_{a,l},H_l)\) if and only if
        $(U_{a,l},H_l)$ is a Nash equilibrium.
  \item Loop mode. Both the strategy profile $(U_a,H)$ and the value function $\dim\operatorname{Ker}\Omega$ evolve on
        a finite cycle, 
        if and only if either player repeats an action after an even number of
        epochs, i.e.\ \(U_{a,i}=U_{a,j}\) or \(H_i=H_j\), whose minimal period divides \((j-i)\).
\end{enumerate}
\end{theorem}
\begin{proof}
(i) (Sufficiency) Assume that the profile $(U_{a,l},H_l)$ is a Nash equilibrium, i.e.
$   U_{a,l}\in BRi^{a}(H_l), \,
   H_l\in BRi^{d}(U_{a,l})$, $i=1 \text{ or }2$.
Because each player is already playing a best response, by
Assumption 1, both players prefer to keep the action unchanged, and inductively in every epoch
$i\ge l$. (Necessity) Conversely, suppose \(\forall i\ge l\)  such that \(U_{a,i}=U_{a,l}\), \(H_i=H_l\).
If $(U_{a,l},H_l)$ were not a Nash equilibrium, at least one
player would have a unilateral deviation in epoch~$l+1$, contradicting that the 
strategy profile $\{U_{a,i},H_i\}$ remains unchanged. Hence $(U_{a,l},H_l)$ must be an equilibrium.

(ii) (Sufficiency) 
Assume that the attacker repeats an action, i.e.\ \(U_{a,i}=U_{a,j}\) with \(j-i\) even.
Then $N+B_a' U_{a,j}=N+B_a' U_{a,i}$,
so the best response set of the defender in epoch \(j+1\) is identical to
that in epoch \(i+1\).
By Assumption~2, the defender
therefore chooses the same action in epoch \(j+1\) as in epoch \(i+1\).
Repeating the argument inductively, we find that each epoch reproduces
the action taken \((j-i)\) periods earlier.
Hence the entire strategy profile and the value function evolve on a
finite cycle whose minimal period divides \(j-i\).
The proof is similar when the defender repeats an action, i.e., \(H_i = H_j\). (Necessity) Since the strategy profile $(U_a,H)$ evolves on a finite cycle, both players repeat the action after the period of the cycle. This completes the proof.
\end{proof}

\subsection{Equilibrium characterization}
Although Theorems \ref{the1} and \ref{the2} provide a unified set of necessary and sufficient conditions for equilibrium, the criterion couples the two players’ strategies, which makes computing the equilibrium directly from the theorem impractical. Therefore, we derive the following necessary condition for the Nash equilibrium. 
\begin{theorem}[Necessary condition for NE]\label{the_necessary}
For any Nash equilibrium $(U_a^*,H^*)$ when considering either one-step or two-step optimality, there must be 
\begin{equation}\label{eq_necessary}
H^* \in \arg \min _{H \in \mathbb{R}^{m\times n}} \operatorname{dim\, }\mathcal{V}^*(H),\,U_a^* \in \mathcal{F}(\mathcal{V}^*(H^*)).
\end{equation}
\end{theorem}
\begin{proof}
According to Theorem \ref{the1} and \ref{the2}, for Nash equilibrium $(U_a^*,H^*)$, there is $\min _{{H}\in \mathbb{R}^{m\times n}}$ $\operatorname{dim\, Ker}\Omega \,(U_a^{*},\!H)=\operatorname{dim}\mathcal{V}^*(H^{*}).$ By further considering Lemma \ref{le4}, we obtain 
\begin{equation}\label{eq_nece}
\min _{H \in \mathbb{R}^{m\times n}} \operatorname{dim\, }\mathcal{V}^*(H)\geq \min _{H \in \mathbb{R}^{m\times n}} \operatorname{dim\, Ker}\Omega(U_a^*,H)=\operatorname{dim}\mathcal{V}^*(H^{*}).
\end{equation}
Since $\operatorname{dim}\mathcal{V}^*(H^{*})\geq \min _{H \in \mathbb{R}^{m\times n}} \operatorname{dim\, }\mathcal{V}^*(H)$, the `$\geq$' symbol in \eqref{eq_nece} in fact holds with equality. Thus $H^* \in \arg \min _{H \in \mathbb{R}^{m\times n}} \operatorname{dim\, }\mathcal{V}^*(H)$. Plus, when considering either one-step or two-step optimality, there must be $U_a^* \in BR1_a(H^*)=\mathcal{F}(\mathcal{V}^*(H^*))$. The proof is completed.
\end{proof}
\begin{remark}
Since $ \min _{H \in \mathbb{R}^{m\times n}} \operatorname{dim\, }\mathcal{V}^*(H)$ in \eqref{eq_necessary} has no relevant to $U_a^*$, we can compute $H^*$ first and then choose $U_a^*$ which belongs to $\mathcal{F}(\mathcal{V}^*(H^*))$. Next, we test whether $(U_a^*,H^*)$ satisfies
$\min_{H\in\mathbb{R}^{m\times n}}
\dim\operatorname{Ker\, }\Omega(U_a^{*},H)=\dim\mathcal{V}^{*}(H^{*})$.
If so, $(U_a^*,H^*)$ is a Nash equilibrium; otherwise, we select another
pair $(U_a^*,H^*)$ that satisfies~\eqref{eq_necessary} and repeat the above procedure
until an equilibrium is found.
\end{remark}
In fact, sometimes Nash equilibrium is disadvantageous to the attacker because the value function remains at a low value that the attacker cannot change. In this case, the attacker can break the equilibrium by forsaking the current gain. Here follows a corollary. 

\begin{corollary}[Attacker-dominated non-equilibrium]\label{the4}
Assume $\min _{H \in \mathbb{R}^{m\times n}}\operatorname{dim}\mathcal{V}^*(H)>0$, which guarantees the existence of at least one attacker strategy $U_a^*$ such that
\begin{equation}\label{eq_nonequli}
U_a^* \notin \mathcal{F}(\mathcal{V}^*(H')), \text{where } H' \in \arg \min _{H \in \mathbb{R}^{m\times n}}\operatorname{dim\, }\mathcal{V}^*(H). 
\end{equation}
Then, for every $H \in \mathbb{R}^{m\times n}$, the profile $(U_a^*,H)$ is not a Nash equilibrium.
\end{corollary}
\begin{proof}
When $\min_{H}\operatorname{dim}\mathcal{V}^*(H)=0$, there is $\mathcal{V}^*(H')=\{0\}$. By the definition of 
\(\mathcal{F}\) we have \(\mathcal{F}(\{0\})=\mathbb{R}^{k\times n}\), so no attacker
strategy can satisfy \eqref{eq_nonequli}.  
When $\min_{H}\operatorname{dim}\mathcal{V}^*(H)>0$, \(\mathcal{V}^*(H')\) is a non-trivial proper subspace of
\(\mathbb{R}^{m}\), implying
\(\mathcal{F}(\mathcal{V}^*(H'))\subsetneq\mathbb{R}\).
Hence \(\mathbb{R}\setminus\mathcal{F}(\mathcal{V}^*(H'))\neq\varnothing\);
choose any
\(U_a^*\) in this complement.  
Such a \(U_a^*\) satisfies \eqref{eq_nonequli} and therefore violates the
necessary condition \eqref{eq_necessary}.
By Theorem~\ref{the_necessary}, $(U_a^*,H)$ fails to satisfy the necessary condition for a Nash equilibrium for all $H\in\mathbb{R}^{m\times n}$. Hence $(U_a^*,H)$ cannot be a Nash equilibrium, which completes the proof.
\end{proof}
By Corollary~\ref{the4}, the attacker can break any candidate equilibrium by
selecting
\(
U_a^*\notin\mathcal{F}\bigl(\mathcal{V}^*(H')\bigr),
\)
where \(\mathcal{V}^*(H')\) is the minimal controlled invariant subspace.
Although this choice forgoes the one–period optimum, it
forces the defender to deviate in the next epoch.  
That deviation enlarges the invariant subspace and raises the attainable
value function.  
The attacker can then switch to an action compatible with the new
subspace, earning a strictly higher utility from the second epoch on and maintaining that advantage thereafter. Section 5 presents an example that illustrates this case.

\subsection{Summary}
To conclude, this chapter analyze the one-step optimality and two-step optimality perspectives and three key insights of the observability-adversarial game can be summarized:

(1) In Subsections 4.1–4.2, Theorems~\ref{the1} and~\ref{the2} show that both the one-step and two-step formulations share a unified necessary-and-sufficient test for a Nash equilibrium, turning different planning horizons into a single easy check. Although the two-step best response sets lie strictly inside the one-step best response sets, the equilibrium reached under either horizon is identical.

(2) In subsection 4.3, as long as the defender selects an \(H^{*}\) that satisfies the defender-dominated Nash equilibrium in Theorem \ref{th3}, the controlled invariant subspace can be collapsed to zero, i.e.\ \(\mathcal V^{*}(H^{*})=\{0\}\). According to Theorem~\ref{cor7}, if in some round the profile \((U_{a},H)\) is a Nash equilibrium, then all subsequent rounds are locked at the same strategy pair (lock mode) and the value function remains fixed, leaving the attacker no further leverage to decrease observability. Hence, the defender not only terminates strategy evolution but also keeps system observability permanently at the most favorable level for the defender.

(3) In subsection 4.4, Theorem~\ref{the_necessary} states the necessary condition for Nash equilibrium.
The attacker can purposely choose a \(U_a^{*}\) violating this condition, sacrificing the immediate best response and thus invalidating the existing locked equilibrium. Once the equilibrium is broken, the strategy trajectory follows Theorem~\ref{cor7} into the Loop mode: both strategies and the value function oscillate on a finite cycle. 
During such a cycle, the controlled invariant subspace can enlarge and system observability may further degrade, so that by accepting a delayed utility the attacker potentially secures a higher long-term value and creates opportunities for deeper penetration.
\section{Illustrative examples}
In this section, we will illustrate the effectiveness of our main results using six cases. Consider a linear system 
\begin{align*}
\dot z&=Nz+B_a'U_az,\\
\hat{y}&=Hz,    
\end{align*}
with $N={ \left[\begin{array}{lllll}
0.3 & 0 & 0 & 0 & 0 \\
0 & 0.3 & 0 & 0 & 0 \\
0 & 0 & 0.3 & 0 & 0 \\
0 & 0 & 0 & 0.1 & 0 \\
0 & 0 & 0 & 0 & 0.2
\end{array}\right]}, B_a'=\left[\begin{array}{lllll}
0 & 0 & 1 & 0 & 1
\end{array}\right]^{\top}$. Initialize with 
$U_a = \begin{bmatrix}0 & 0 & 0 & 0 & 0\end{bmatrix}$.
The defender then selects 
$H\in\mathbb{R}^{2\times 5}$
to minimize the dimension of the unobservable subspace. Next, the attacker updates 
$U_a\in\mathbb{R}^{1\times 5}$
to maximize that dimension. This sequence of moves repeats and yields a repeated game. Although only the outcomes of 60 epochs are shown in the following figures, the cyclic nature of the results makes these snapshots sufficient to characterize the behavior over an infinite time horizon.

The game settings for Case 1-6 are summarized in Table 2. Cases are divided into three categories (Cases 1-3; Cases 4-5; Case 6), which will be introduced separately below.
\begin{center}
{\small
    \begin{center}
\centerline{\small {\bf Table 2}\ \  ~~Game settings for Cases 1–6}
    \label{tab:my_label}
    {\small
    \begin{tabular*}{\textwidth}{|p{1cm}|p{5cm}|p{5cm}|p{1.8cm}|}
    \hline
    Settings & Attacker action rules & Defender action rules & Strategies \\              
    \hline
    Case 1 & Algorithm 1. & Algorithm 2. & \multirow{3}{=}{One‐step optimality} \\
    \cline{1-3}
    Case 2 & \renewcommand{\arraystretch}{1} 
      $U_a\!=\!\bigl[\begin{array}{lllll}
        0 & 0 & 0 & 0.2 & 0.1
      \end{array}\bigr]$
      instead of 
      $\bigl[\begin{array}{lllll}
        -0.1 & 0 & 0 & 0.1 & 0
      \end{array}\bigr]$
      when 
      $H\!=\!\bigl[\begin{array}{lllll}
        1 & 0 & 0 & 1 & 1 \\
        0 & 1 & 0 & 0 & 0
      \end{array}\bigr],$ 
      and other actions based on Algorithm 1. 
      & Algorithm 2. & \\
    \cline{1-3}
    Case 3 & Algorithm 1. 
      & \renewcommand{\arraystretch}{1} 
      $H\!=\!\bigl[\begin{array}{lllll}
        1 & 0 & 0 & 1 & 1 \\
        0 & 0 & 1 & 0 & 0
      \end{array}\bigr]$
      instead of 
      $\bigl[\begin{array}{lllll}
        1 & 0 & 0 & 1 & 1 \\
        0 & 1 & 0 & 0 & 0
      \end{array}\bigr]$
      when 
      $U_a\!=\!\bigl[\begin{array}{lllll}
        0 & 0 & 0 & 0 & 0
      \end{array}\bigr],$ 
      and other actions follow Algorithm 2. 
      & \\
    \hline
    Case 4 & Before epoch 20: Algorithm 1; after 20 epochs: two‐step optimality. 
      & Before epoch 40: Algorithm 2; after 40 epochs: two‐step optimality. 
      & \multirow{2}{=}{One‐step \& two‐step optimality} \\
    \cline{1-3}
    Case 5 & Before epoch 40: Algorithm 1; after 40 epochs: two‐step optimality.
      & Before epoch 20: Algorithm 2; after 20 epochs: two‐step optimality.
      & \\
    \hline
    Case 6 & One‐step optimality, only with $U_a \notin \mathcal{F}(\mathcal{V}^*)$ in epoch 40.
      & One‐step optimality.
      & One‐step optimality \& attacker not greedy \\
    \hline
    \end{tabular*}
    }
    \end{center}
}
\end{center}

\textbf{Cases 1-3:} let players have different choices in the one-step optimal best response sets. In Case 1, the attacker chooses $U_a\!=\!\left[\begin{array}{lllll}
\!-0.1 & 0 & 0 & 0.1 & 0
\end{array}\right]$ when $H\!=\!\left[\begin{array}{lllll}1 & 0 & 0 & 1 & 1 \\
0 & 1 & 0 & 0 & 0
\end{array}\right]$ according to Algorithm 1, which uses pseudo inverse to get $U_a$ with the minimum modulus length; the defender chooses $H=\left[\begin{array}{lllll}
1 & 0 & 0 & 1 & 1 \\
0 & 1 & 0 & 0 & 0
\end{array}\right]$ when $U_a\!=\!\left[\begin{array}{lllll}
0 & 0 & 0 & 0 & 0
\end{array}\right]$ according to Algorithm 2. While in Cases 2-3, players have different choices compared to Case 1. 
Figure 5 shows the evolution of unobservable subspace dimension for Cases 1-3, which are different for different cases. The specific actions of players $U_a$ and $H$ for cases 1-3 are summarized in Table 3. In cases 1-2, the actions $U_a$ and $H$ evolve in different loops, which illustrate the results of Theorem \ref{cor7}. In Case 3, the actions $U_a^*$ and $H^*$ keep unchanged, which is the candidate equilibrium. We then verify that
 $ \min_{H\in\mathbb{R}^{m\times n}}\operatorname{dim\, Ker}\Omega(U_a^*,H)
  \;=\;\dim\mathcal{V}^*(H^*)$,
thereby confirming the validity of Theorem~\ref{the1} and Theorem \ref{cor7}.

\begin{center}
  \centerline
  {\includegraphics[scale=0.6]{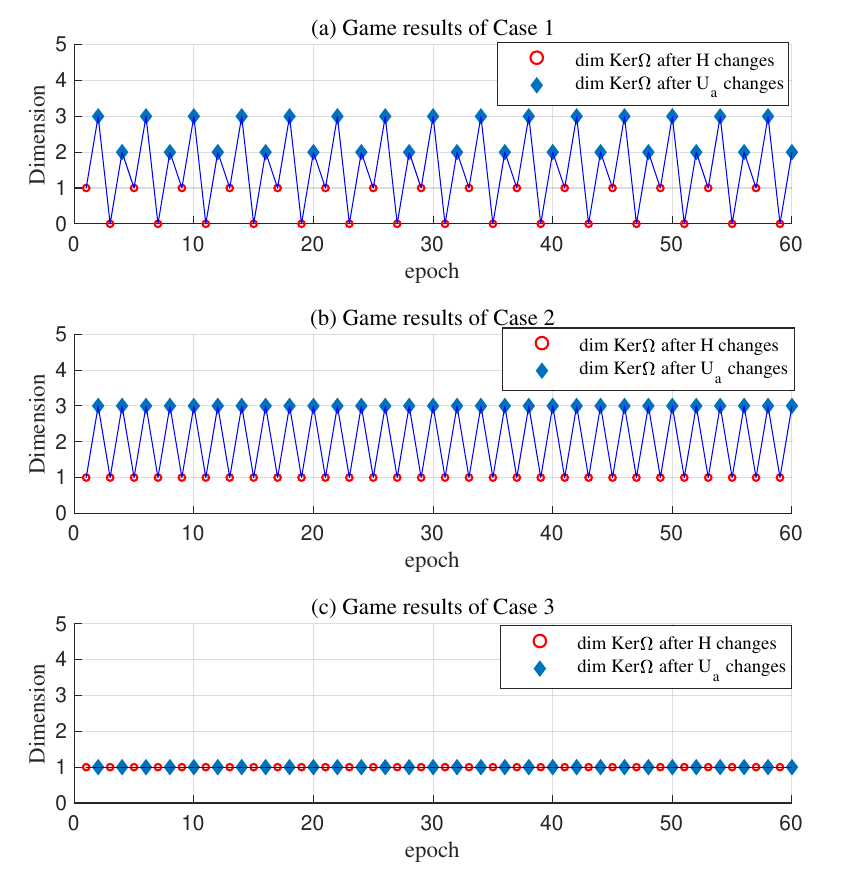}}
\centering{\small {\bf Figure 5}\ \ Evolution of $\operatorname{dim\, Ker}\Omega$ for Cases 1-3. (a) Case 1: the attacker uses Algorithm 1 and the defender uses Algorithm 2; (b) Case 2: the attacker changes strategy preference and other actions match Case 1; (c) Case 3: the defender changes strategy preference and other actions match Case 1. }
\end{center}

\begin{table*}[htbp] 
\centering 
\centerline{\small {\bf Table 3}\ \  Actions evolution of players for Cases 1-3.} 
\renewcommand{\arraystretch}{1.5}
\footnotesize
\begin{tabular}{|p{0.9cm}|p{4.1cm}|p{4cm}|p{3.8cm}|}
\hline
Epochs $(n\in N^*)$ & Case 1 & Case 2 & Case 3 \\
\hline
$n$ & 
\renewcommand{\arraystretch}{1}
$H_n=\begin{bmatrix}
1 & 0 & 0 & 1 & 1\\
0 & 1 & 0 & 0 & 0
\end{bmatrix}$ &
\renewcommand{\arraystretch}{1}
$H_n=\begin{bmatrix}
1 & 0 & 0 & 1 & 1\\
0 & 1 & 0 & 0 & 0
\end{bmatrix}$ &
\renewcommand{\arraystretch}{1}
$H_n=\begin{bmatrix}
1 & 0 & 0 & 1 & 1\\
0 & 0 & 1 & 0 & 0
\end{bmatrix}$ \\
2$n$ &
\renewcommand{\arraystretch}{1}
$U_{a2n}\!=\!\begin{bmatrix}
-0.1 & 0 & 0 & 0.1 & 0
\end{bmatrix}$ &
\renewcommand{\arraystretch}{1}
$U_{a2n}=\begin{bmatrix}0 & 0 & 0 & 0.2 & 0.1\end{bmatrix}$&
\renewcommand{\arraystretch}{1}
$U_{a2n}=\begin{bmatrix}
0 & 0 & 0 & 0 & 0
\end{bmatrix}$ \\
3$n$ &
\renewcommand{\arraystretch}{1}
$H_{3n}=\begin{bmatrix}
0 & 0 & -1 & 0 & 1\\
0 & 1 & 0 & 0 & 0
\end{bmatrix}$  &
\renewcommand{\arraystretch}{1}
$H_{3n}=\begin{bmatrix}
0 & 0 & 1 & 0 & 0\\
1 & 0 & 0 & 0 & 0
\end{bmatrix}$  &
\renewcommand{\arraystretch}{1}
$H_{3n}=\begin{bmatrix}
1 & 0 & 0 & 1 & 1\\
0 & 0 & 1 & 0 & 0
\end{bmatrix}$  \\
4$n$ &
\renewcommand{\arraystretch}{1}
$U_{a4n}=\begin{bmatrix}
0 & 0 & 0 & 0 & 0
\end{bmatrix}$ &
\renewcommand{\arraystretch}{1}
$U_{a4n}=\begin{bmatrix}
0 & 0 & 0 & 0 & 0
\end{bmatrix}$ &
\renewcommand{\arraystretch}{1}
$U_{a4n}=\begin{bmatrix}
0 & 0 & 0 & 0 & 0
\end{bmatrix}$ \\
\hline
\end{tabular}
\end{table*}

\textbf{Cases 4-5:} We allow two players to apply two-step optimality in different epochs to compare their impact on game result. Before epoch 20, let two players consider one-step optimality and choose actions according to Algorithm 1 and 2. In Case 4, let the attacker consider two-step optimality after 20 epochs and both players consider two-step optimality after 40 epochs. In Case 5, let the defender consider two-step optimality after 20 epochs and other settings are the same as Case 4.
\begin{center}
  \begin{minipage}[t]{0.47\linewidth}
    \centering
    \includegraphics[scale=0.42]{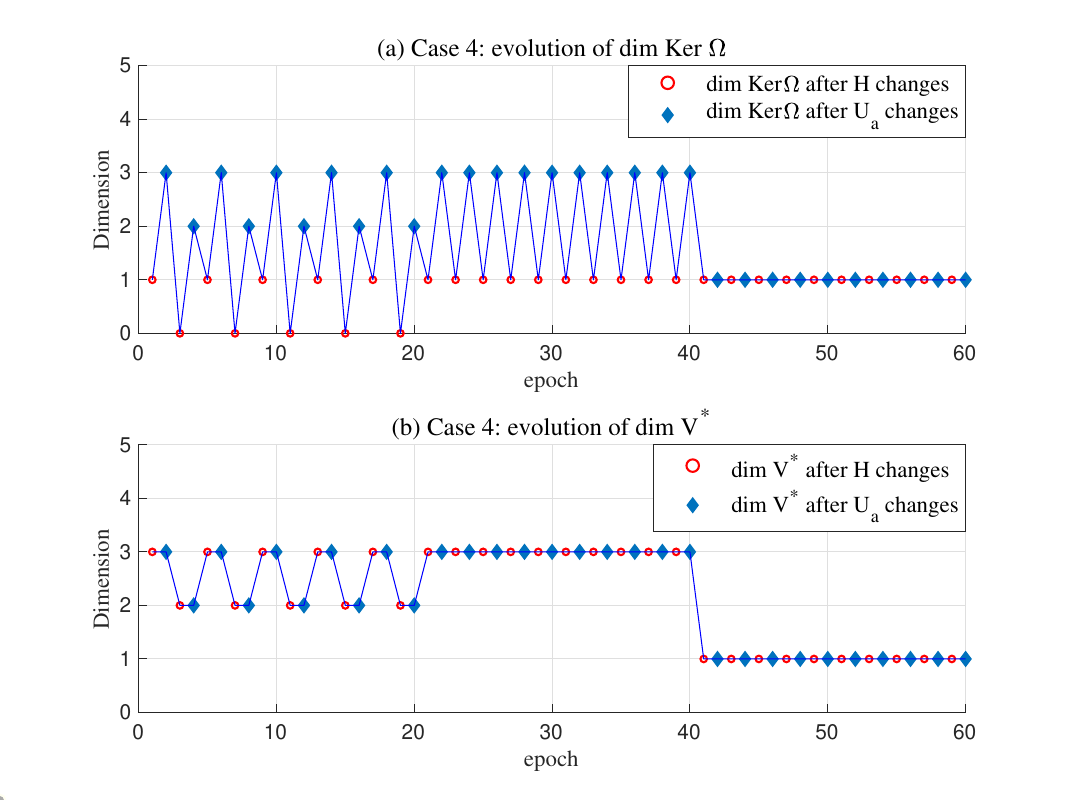}
    \vskip3mm
    {\small {\bf Figure 6}\ \ Results of Case 4: the attacker considers two-step optimality after 20 epochs. (a) Evolution of $\operatorname{dim\, Ker}\Omega$; (b) Evolution of $\dim\mathcal{V}^*$.}
  \end{minipage}%
  \hfill
  \begin{minipage}[t]{0.47\linewidth}
    \centering
    \includegraphics[scale=0.42]{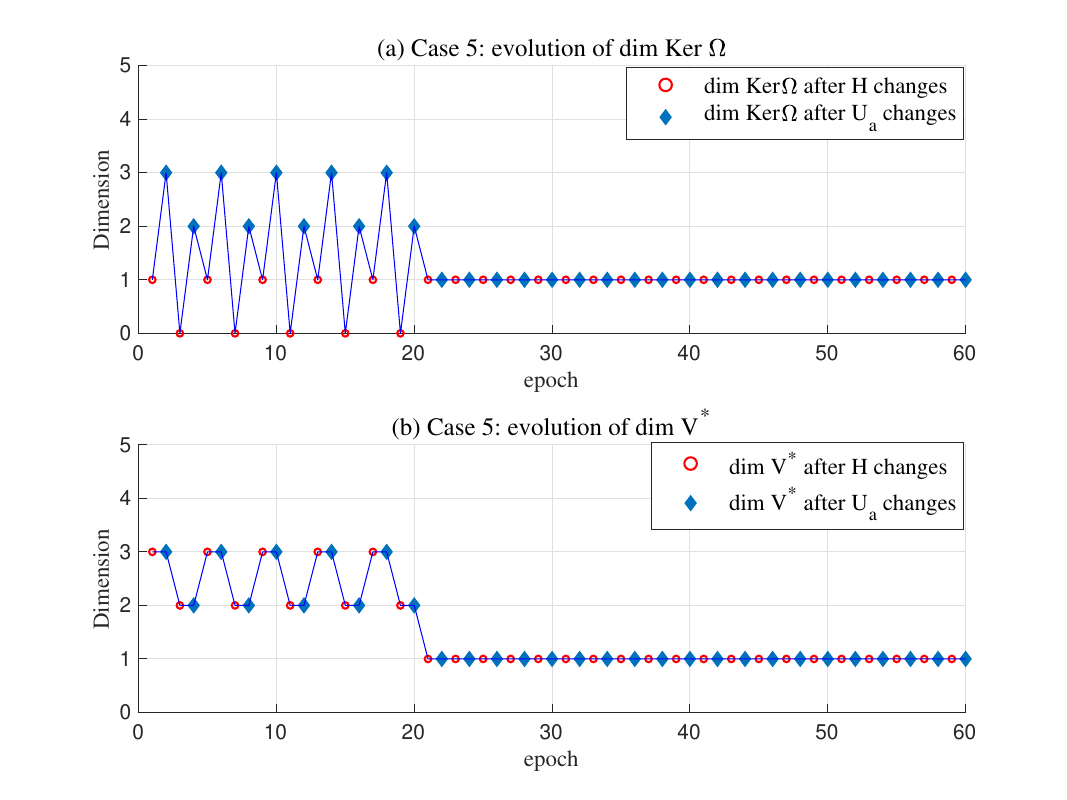}
    \vskip3mm
    {\small {\bf Figure 7}\ \ Results of Case 5: the defender considers two-step optimality after 20 epochs. (a) Evolution of $\operatorname{dim\, Ker}\Omega$; (b) Evolution of $\dim\mathcal{V}^*$.}
  \end{minipage}
\end{center}

Figures 6–7 depict Cases 4–5. Up to epoch 20 both cases follow the
trajectory of Case 1.  At epoch 20 the attacker (Case 4) or the defender (Case 5) adopts a two-step best response, producing different outcomes.
In Case 4, \(\dim\mathcal{V}^{*}\) jumps to 3, exceeding  
  \(\min_{H}\operatorname{dim\, Ker}\Omega(U_a,H)=1\); no Nash equilibrium exists and the
  state cycles until epoch 40. 
In Case 5, \(\dim\mathcal{V}^{*}\) drops to 1 and equals  
  \(\min_{H}\operatorname{dim\, Ker}\Omega(U_a,H)\); the game reaches a Nash equilibrium
  immediately.
After epoch 40 both cases satisfy  
\(\dim\mathcal{V}^{*}=\min_{H}\operatorname{dim\, Ker}\Omega(U_a,H)=1\) and remain at
equilibrium.  These results corroborate the necessary and sufficient condition of NE in Theorem~\ref{the2}. 

Furthermore, in Case 4 the attacker’s two-step move raises the long-run average of $\operatorname{dim\, Ker}\Omega$ from 1.5 to roughly 2.0, whereas in Case 5 the defender’s two-step response lowers it from 1.5 to roughly 1.0. Thus, by selecting an action from its two-step optimal set, either player can shift the long-run value function in its own favor.

In the above five cases, both players are greedy, i.e., their actions are either one-step or two-step optimal. Next, consider a case when the attacker is not greedy.

\textbf{Case 6:} let both players consider one-step optimality and choose actions to achieve equilibrium before epoch 40. In epoch 40, when $H=\left[\begin{array}{lllll}
1 & 0 & 0 & 1 & 1 \\
0 & 0 & 1 & 0 & 0
\end{array}\right]$, the attacker chooses $U_a=\left[\begin{array}{lllll}
0 & 1 & 0 & 0 & 0
\end{array}\right] \notin BR1^a$. 
\begin{center}
  \centerline
  {\includegraphics[scale=0.42]{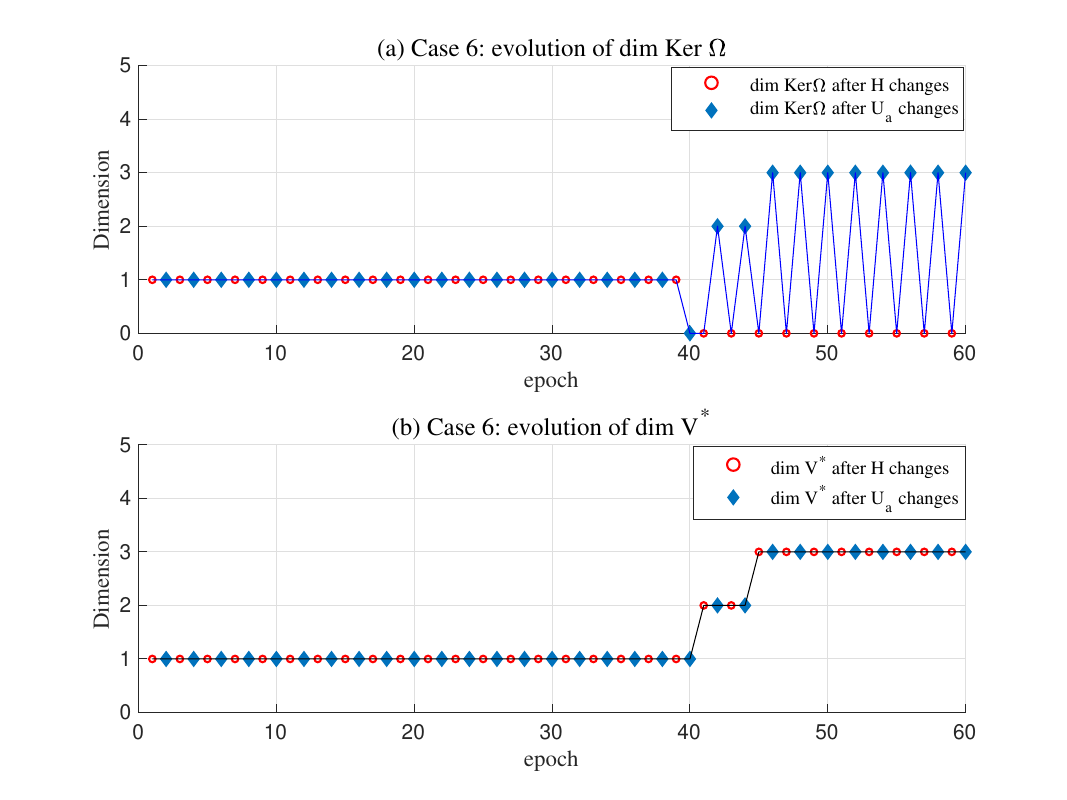}}\vskip3mm
\centering{\small {\bf Figure 8}\ \ Results of Case 6: in epoch 40, the attacker chooses $U_a \notin  BR1^a$. (a) Evolution of $\operatorname{dim\, Ker}\Omega$; (b) Evolution of $\dim\mathcal{V}^*$.}
\end{center}

The game result of case 6 is shown in Figure 8. Although the value function in epoch 40 reduces, $\operatorname{dim\,}\mathcal{V}^*$ in epoch 41 increases, which increases the value function in epoch 42. Thus the equilibrium is broken and the game reaches a new mode: the value function oscillates between 1 and 3, which illustrates the result of Corollary \ref{the4}. From the Nash equilibrium to the oscillation mode, the average value function increases from 1 to 1.5, which is more beneficial for the attacker. It shows that by giving up the benefit of the current epoch, the attacker can break the Nash equilibrium and make the game result more beneficial for itself. 
\section{Conclusion}
In this paper, we formulate the contest between an attacker and a defender over system observability as an infinitely repeated game whose value function equals the dimension of the unobservable subspace. Derivations and algorithms to maximize or minimize the unobservable subspace dimension are given. Despite the resulting best response sets being set-valued, we established a unified necessary-and-sufficient condition for Nash equilibrium. 
The long-term outcome of the game depends on whether the two players satisfy the Nash equilibrium conditions. 
If a defender-dominated Nash equilibrium exists and the defender chooses the corresponding strategy, the controlled invariant subspace collapses immediately to $\{0\}$ and the game enters the lock mode. In this case, the attacker can no longer reduce observability and system security remains permanently fixed at the most favorable level of the defender.  
Under the more general non-degenerate condition, the game admits two possible outcomes: lock mode and loop mode.
Furthermore, we provide a necessary condition for the Nash equilibrium in which the strategies of two players are uncoupled. By deliberately adopting a strategy that violates this condition, the attacker can break the equilibrium. Although the attacker sacrifices an immediate utility, it may achieve a higher long-term value. Finally, numerical case studies confirm these insights.

In the future, we can study the problem in more complex scenarios. Factors such as system stability, energy consumption and unobservable subspace can be combined to formulate a more comprehensive value function. Moreover, we can consider the game with incomplete information, for example with unknown system matrices, how does the attacker design strategies to change the observability of the system, and the strategy evolution of both sides.
\renewcommand\bibname{References}


\begin{thebibliography}{100}

\bibitem{niu2023innovation}
Niu M, Wen G, Lv Y, et al., Innovation-based stealthy attack against distributed state estimation over sensor networks, {\it Automatica}, 2023, {\bf 152}: 110962.

\bibitem{soltan2018react}
Soltan S, Yannakakis M, and Zussman G, REACT to cyber attacks on power grids, {\it IEEE Trans. Netw. Sci. Eng.}, 2018, {\bf 6}(3): 459–473.

\bibitem{bernard2022observer}
Bernard P, Andrieu V, and Astolfi D, Observer design for continuous-time dynamical systems, {\it Annu. Rev. Control}, 2022, {\bf 53}: 224–248.

\bibitem{showkatbakhsh2020securing}
Showkatbakhsh M, Shoukry Y, Diggavi S N, et al., Securing state reconstruction under sensor and actuator attacks: Theory and design, {\it Automatica}, 2020, {\bf 116}: 108920.

\bibitem{liu2024performance}
Liu Q, Wang J, Ni Y, et al., Performance analysis for cyber–physical systems under two types of stealthy deception attacks, {\it Automatica}, 2024, {\bf 160}: 111446.

\bibitem{chong2015obser}
Chong M S, Wakaiki M, and Hespanha J P, Observability of linear systems under adversarial attacks, {\it 2015 American Control Conference}, 2015, pp. 2439–2444.

\bibitem{fawzi2014secure}
Fawzi H, Tabuada P, and Diggavi S, Secure estimation and control for cyber-physical systems under adversarial attacks, {\it IEEE Trans. Automatic Control}, 2014, {\bf 59}(6): 1454–1467.

\bibitem{shoukry2015event}
Shoukry Y and Tabuada P, Event-triggered state observers for sparse sensor noise/attacks, {\it IEEE Trans. Automatic Control}, 2015, {\bf 61}(8): 2079–2091.

\bibitem{mitra2018distributed}
Mitra A and Sundaram S, Distributed observers for LTI systems, {\it IEEE Trans. Automatic Control}, 2018, {\bf 63}(11): 3689–3704.

\bibitem{mitra2019byzantine}
Mitra A and Sundaram S, Byzantine-resilient distributed observers for LTI systems, {\it Automatica}, 2019, {\bf 108}: 108487.

\bibitem{zhang2023observability}
Zhang Y, Xia Y, and Liu K, Observability robustness under sensor failures: A computational perspective, {\it IEEE Trans. Automatic Control}, 2023, 68(12): 8279-8286.

\bibitem{kim2014subspace}
Kim J, Tong L, and Thomas R J, Subspace methods for data attack on state estimation: A data driven approach, {\it IEEE Trans. Signal Process.}, 2014, {\bf 63}(5): 1102–1114.

\bibitem{zhao2019sparse}
Zhao Z, Li Y, Yang Y, et al., Sparse undetectable sensor attacks against cyber-physical systems: A subspace approach, {\it IEEE Trans. Circuits Syst. II: Express Briefs}, 2019, {\bf 67}(11): 2517–2521.

\bibitem{maccarone2018uncovering}
Maccarone L T, D’Angelo C J, and Cole D G, Uncovering cyber-threats to nuclear system sensing and observability, {\it Nuclear Eng. Des.}, 2018, {\bf 331}: 204–210.

\bibitem{zhoupanpan}
Zhou P and Chen B M, Distributed Optimal Solutions for Multiagent Pursuit-Evasion Games for Capture and Formation Control, {\it IEEE Trans. Ind. Electron.}, 2024, {\bf 71}(5): 5224–5234.

\bibitem{zhang2025}
Zhang Y, Lian B, Lewis F L, Distributed global nash equilibrium of interactive adversarial graphical games, {\it Journal of Systems Science and Complexity}, 2025, {\bf 38}(2): 613-632.

\bibitem{horak2019optimizing}
Horák K, Bosanský B, Tomášek P, et al., Optimizing honeypot strategies against dynamic lateral movement using partially observable stochastic games, {\it Computers \& Security}, 2019, {\bf 87}: 101579.

\bibitem{zheng2022stackelberg}
Zheng W, Jung T, and Lin H, The Stackelberg equilibrium for one-sided zero-sum partially observable stochastic games, {\it Automatica}, 2022, {\bf 140}: 110231.

\bibitem{maccarone2020game}
Maccarone L T and Cole D G, A game-theoretic approach for defending cyber-physical systems from observability attacks, {\it ASCE-ASME J. Risk Uncertainty Eng. Syst., Part B: Mec. Eng.}, 2020, {\bf 6}(2): 021004.

\bibitem{nguyen2019deception}
Nguyen T H, Wang Y, Sinha A, et al., Deception in finitely repeated security games, {\it Proceedings of the AAAI Conference on Artificial Intelligence}, 2019, 33(01): 2133-2140.

\bibitem{balaji2019design}
Balaji S, Julie E G, Robinson Y H, et al., Design of a security-aware routing scheme in mobile ad-hoc network using repeated game model, {\it Computer Standards \& Interfaces}, 2019, {\bf 66}: 103358.

\bibitem{aziz2020resilience}
Aziz F M, Li L, Shamma J S, et al., Resilience of LTE eNode B against smart jammer in infinite-horizon asymmetric repeated zero-sum game, {\it Physical Communication}, 2020, {\bf 39}: 100989.

\bibitem{xie2025}
Xie K, Lu M, Deng F, et al., Data-Driven Dynamic Output Feedback Nash Strategy for Multi-Player Non-Zero-Sum Games, {\it Journal of Systems Science and Complexity}, 2025, {\bf 38}(2): 597-612.

\bibitem{rizvi2025}
Anwar J, Rizvi S A A, Lin Z, Output Feedback Q-Learning for a Non-Zero-Sum Game Problem in Building HVAC Control, {\it Journal of Systems Science and Complexity}, 2025, {\bf 38}(2): 739-755.

\bibitem{xu2023}
Xu Y, Hu X, Liu Z, et al., A Game Approach for Defending System Security from an Attacker, {\it IFAC-PapersOnLine}, 2023, {\bf 56}(2): 1710–1715.

\bibitem{Teixeira2012}
Teixeira A, Shames I, Sandberg H, et al., Revealing stealthy attacks in control systems, {\it 2012 50th Annual Allerton Conference on Communication, Control, and Computing}, 2012, pp. 1806–1813.

\bibitem{basile1969}
Basile G and Marro G, Controlled and conditioned invariant subspaces in linear system theory, {\it J. Optim. Theory Appl.}, 1969, {\bf 3}: 306–315.

\bibitem{chen1984linear}
Chen C T, {\it Linear system theory and design}, Saunders College Publishing, 1984.
\end{thebibliography}
\end{document}